\documentclass[11pt]{amsart}
\usepackage{cite}
\usepackage{wrapfig, lipsum,booktabs}
\usepackage{graphicx}
\usepackage{amssymb}
\usepackage{epstopdf}
\usepackage{verbatim}
\usepackage{bm}
\usepackage{multicol}
\usepackage{multirow}
\usepackage{subfigure}
\usepackage{fancyhdr} 
\usepackage{float}
\usepackage{stmaryrd}
\usepackage{color}
\usepackage{soul}
\usepackage{tikz}
\usetikzlibrary{patterns}
\usepackage{pgfplots}
\usepackage{cancel}

\newtheorem{remark}{Remark}[section]
\usepackage{multirow}
\usepackage{amsmath,amssymb,eucal}
\usepackage{graphicx,subfigure,epsfig}
\usepackage{psfrag}
\usepackage{url}
\usepackage{todonotes}
\usepackage[top=1in, bottom=1.in, left=1in, right=1in]{geometry}

\newcommand{\bld}[1]{\hbox{\boldmath$#1$}}    
\newcommand{\Eh}{\mathcal{E}_h}
\newcommand{\jmp}[1]{[\![#1 ]\!]}
\def\avg#1{\{\!\!\{{#1}\}\!\!\}} 
\newcommand{\Oh}{\mathcal{T}_h}
\newcommand{\pol}{\mathbb{P}}
\newcommand{\Vhat}[1]{\widehat{\bld V}_{\!\!h}^{#1}}
\newcommand{\What}[1]{\widehat{W}_{\!h}^{#1}}
\newcommand{\Vh}[1]{{\bld V}_{\!\!h}^{#1}}
\newcommand{\Wh}[1]{{W}_{\!h}^{#1}}

\begin{document}
\title[POD-(H)DG for incompressible flows]{POD-(H)DG method for incompressible flow simulations}
\author{Guosheng Fu}
\address{Department of Applied and Computational Mathematics and 
Statistics, University of Notre Dame, USA.}
\email{gfu@nd.edu}
\author{Zhu Wang}
\address{Department of Mathematics, University of South Carolina, Columbia, SC 29208, USA.}
\email{wangzhu@math.sc.edu}
 \thanks{Zhu Wang gratefully acknowledges the partial support of this work
 from U.S. National Science Foundation through grant DMS-1913073.}

\begin{abstract}
  We present a reduced order method (ROM) based on proper orthogonal decomposition
(POD) for the viscous Burgers' equation and the incompressible Navier-Stokes equations
discretized using an implicit-explicit hybrid discontinuous
Galerkin/discoutinuous Galerkin (IMEX HDG/DG) scheme. 
A novel closure model, which can be easily computed offline,
is introduced. Numerical results are presented to test the proposed POD
model and the closure model. 
\end{abstract}
\maketitle

\section{Introduction}
\label{sec:intro}

Reduced order modeling has been widely used in flow control and optimization problems to alleviate the huge computational cost needed in many-query solutions of the large-scale dynamical systems associated to these problems \cite{gunzburger2003perspectives,antoulas2005approximation,quarteroni2015reduced,hesthaven2016certified}. To achieve the computationally high efficiency, model reduction methods construct from data a numerical surrogate model with the dimension greatly reduced from the original system. To build such a low-dimensional model, one can use non-intrusive approaches such as operator learning  \cite{peherstorfer2016data,antoulas2020interpolatory}, or intrusive approaches such as projection-based methods \cite{benner2015survey}. The method to be used in this paper falls into the second category. In particular, we consider the proper orthogonal decomposition (POD) method - one of the most popular snapshot-based model reduction techniques. The general POD model reduction methodology splits the overall calculation into offline and online stages. At the offline stage, a handful of reduced basis vectors are determined and a low-dimensional, reduced order model (ROM) is constructed by learning algorithms or by projecting equations to the space spanned by the reduced basis. At the online stage, the ROM is used alternative to the original system for simulations that can be finished in short time or even real time. When the system contains non-polynomial nonlinearities, hyper-reduction has to be used in order to guarantee the online computational complexity to be independent of the dimension of the original system \cite{chaturantabut2010nonlinear,carlberg2013gnat}.

The ROM can be discretized by any conventional numerical method. In particular, when continuous Galerkin finite elements are used, each nodal value will be shared by several elements. If an interpolation type of hyper-reduction methods is applied, such as discrete empirical interpolation method (DEIM) or its variants \cite{chaturantabut2010nonlinear}, although the nonlinear functions need only to evaluate at few selected points, many elements that share these nodes have to be looped. This would cause expensive online computations. Thus, the finite element with interpolated coefficients method was developed in \cite{wang2015nonlinear}, in which the nonlinear functions in the ROM are replaced with their finite element interpolants so that the DEIM can be applied directly on the finite element coefficients. However, if a discontinuous Galerkin (DG) method is applied, there is no such issue thanks to the local nature of the DG method. 
There has been several work that uses POD in the context of DG. In \cite{shen2019hdg} hybridizable discontinuous Galerkin (HDG) POD model has been developed for heat equation. It is shown that highly accurate flux approximation can be recovered in the HDG-POD approximation at a low cost. In \cite{uzunca2017energy}, POD is applied in the context of symmetric interior penalty DG for solving Allen-Cahn equation. For parametric problems, DG has also been applied together with reduced basis method for elliptic problems \cite{antonietti2016discontinuous} and with empirical quadrature procedure for nonlinear conservation laws in \cite{yano2019discontinuous}.

As a first step for investigating reduced order modeling on flow control
and optimization applications, we focus on the computational fluid dynamics
of incompressible fluid flows in this work. When the POD approximation is
sought for such problems, there are two common ways  to deal with the
incompressibility constraint. One only keeps velocity in the reduced
system, which is based on the argument that the POD basis is weakly
divergence-free since it is a combination of snapshots and snapshots are
weakly divergence-free, thus the pressure term would vanish after
projection; the other keeps both velocity and pressure in the reduced
system, because either the application at hand needs pressure information
or numerical methods for computing snapshots may not provide pointwise
divergence free flow fields. Indeed, the discretely divergence-free
property does not hold for many popular discretization of the Navier-Stokes
equations. A new velocity and pressure ROM is proposed in \cite{ballarin2015supremizer} while introducing a supremizer stabilization to fulfill an equivalent inf-sup condition. In \cite{caiazzo2014numerical}, one velocity ROM and two velocity-pressure ROMs are compared that shows the accuracy of snapshots does have a big impact on the performance of velocity ROM. 
Therefore, in this work, we use the divergence-free HDG method developed in
\cite{Lehrenfeld10} for the full order mode (FOM), which ensures the
velocity snapshots are exactly pointwise divergence free. As a consequence,
we can use the velocity ROM since the POD basis generated from these
snapshots would have the same divergence-free property. Furthermore,
because the convective term involves upwinding numerical flux in the FOM
that can not be precomputed offline, we replace it with
an (offline-computable) linear central flux in the ROM. 
We then add a {\it linear} convective stabilization term in the ROM which
mimics the upwinding stabilization in the FOM. 
This yields an efficient implementation while keeping stabilization effects of the numerical flux. However, the introduced stabilization might not be enough for convection-dominated problems, because the jumps across elements of the POD basis are small when snapshots are obtained from high fidelity simulations. Therefore, we include extra dissipation following the closure model developed in \cite{SanIliescu13} to diminish the numerical oscillations. 

The rest of paper is organized as follows. In Section \ref{sec:fom}, the
full order model is presented. In Section \ref{sec:pod}, the POD reduced
order model is derived. Several numerical experiments including the
Burgers' equation, the Navier-Stokes equations and the incompressible Euler
equations are discussed in Section \ref{sec:num}. 
We conclude in Section \ref{sec:conclude} with some future work.

\section{Full order model via IMEX HDG/DG}
\label{sec:fom}

Next, we first describe the FOM that is employed for generating snapshots and provides the benchmark solutions in our numerical experiments. 
\subsection{Notation} 
Let $\Oh$ be a conforming simplicial triangulation of 
the domain $\Omega\subset
\mathbb{R}^d$, $d=1,2,3$.
For any element $K \in\Oh$, we denote by $h_K$ its diameter and by $h:
\Oh\rightarrow \mathbb{R}$ the mesh size function with $h|_K = h_K$. 
The collection of element boundaries is $\partial\Oh:=\{\partial K:\;
K\in\Oh\}$.
Denote by $\Eh$ the set of facets of  $\Oh$ (vertices in 1D, edges in 2D,
faces in 3D), and by $\Eh^i=\Eh\backslash\partial\Omega$
the set of interior facets. 
For any element $K$, denote $\bld n_{K}:\partial K\rightarrow \mathbb{R}^d$
to be the unit outward normal direction on $\partial K$ from the element $K$. 
Let  $\bld n:\partial\Oh\rightarrow\mathbb{R}^d$ be the unit normal direction on the collection of element
boundaries $\partial \Oh$ with $\bld n|_{\partial K} = \bld n_K$.

We collect the following set of finite element spaces:
\begin{subequations}
\label{space}
\begin{align}
\label{space-1}
\Vh{k} : =&\; \{\bld v\in 
  H(\mathrm{div},\Omega):\;\;
  \bld v|_{K}\in [\pol^{k}(K)]^d,\;\; \forall K\in \Oh,\;
(\bld v\cdot\bld n)|_{\partial\Omega} = 0\},\\
\label{space-2}
  \Vhat{k} := &\;\{\widehat{\bld v}\in [L^2(\Eh)]^d:
    \;\;\widehat{\bld v}|_F\in  [\pol^k(F)]^d, \;
    (\widehat{\bld v}\cdot\bld n)|_F = 0, \;\;\forall F\in\Eh, \;\;
  \widehat{\bld v}|_{\partial\Omega} = 0\},\\
\label{space-3}
    \Wh{k} : =&\; \{w\in 
  L^2(\Omega):\;\;
  w|_{T}\in \pol^{k}(K),\;\; \forall K\in \Oh\},\\
\label{space-4}
  \What{k} := &\;\{\widehat{w}\in L^2(\Eh):
    \;\;\widehat{w}|_F\in  \pol^k(F), \;
  \widehat{w}|_{\partial\Omega} = 0\},
\end{align}
\end{subequations}
where 
$\pol^k$ is the space of polynomials up to degree $k\ge 0$.
In 1D, $\pol^k(F)$ is simply point evaluation for the vertex 
$F$.
Note that functions in $\Vhat{k}$ and $\What{k}$ are defined only on the mesh skeleton $\Eh$. 


\subsection{The model problems}
Two mathematical models are considered in this work, namely the 1D viscous Burgers'
equation
  \begin{alignat}{2}
  \label{eq:burgers}
  \frac{\partial u}{\partial t}+
  u\cdot\nabla u
  -\nu\triangle u=&\;0 &&\quad \text{ in }\Omega\subset \mathbb{R},
\end{alignat}
and the 2D incompressible Navier-Stokes
equations \eqref{eq:ns}:
\begin{subequations}
  \label{eq:ns}
\begin{align}
  \frac{\partial \bld u}{\partial t}+
  \bld u\cdot\nabla \bld u + \nabla p
  -\nu\triangle \bld u=&\;0 \quad \text{ in }\Omega\subset\mathbb{R}^2,\\
  \nabla \cdot\bld u 
  =&\;0 \quad \text{ in }\Omega,
\end{align}
\end{subequations}
where $\nu> 0$ is a positive viscosity parameter. It becomes the incompressible Euler equations when $\nu=0$. 
For simplicity of presentation, we use homogeneous Dirichlet boundary
conditions for both problems to derive the FOM IMEX HDG/DG schemes. 
Other standard boundary conditions will be applied in the numerical experiments presented in 
Section \ref{sec:num}.

\subsection{The semidiscrete HDG/DG scheme: Burgers' equation}
The semidiscrete HDG/DG scheme for the 1D Burgers' equation
\eqref{eq:burgers}
reads as follows:
Given initial data $u_h(0)\in \Wh{k}$, for all 
$t\in (0,T]$, find $(u_h,\widehat{u}_h)=
(u_h(t), \widehat{u}_h(t))\in \Wh{k}\times \What{0}$ such
that 
\begin{align}
  \label{semi-b}
  \mathcal{M}_h(\frac{\partial u_h}{\partial t}, v)
  +\mathcal{C}_h^{dg}(u_h,u_h,v)
  +\nu\mathcal{B}_h^{hdg}((u_h, \widehat{u}_h),
  (v, \widehat{v})) = 0, \quad\forall 
  (v,\widehat{v})\in \Wh{k}\times \What{0}.
\end{align}
Here $\mathcal{M}_h(\cdot,\cdot)$ is the mass operator, 
$\mathcal{C}_h^{dg}(\cdot, \cdot,\cdot)$ is the nonlinear (DG) convection operator, 
and $\mathcal{B}_h^{hdg}(\cdot,\cdot)$ is the (HDG) diffusion operator, which are
given as follows:
\begin{subequations}
  \label{ops-b}  
  \begin{alignat}{2}
  \label{mass-b}
  \mathcal{M}_h(u,v)=&\;
  \sum_{K\in\Oh}\int_{K}u\,v\,\mathrm{dx},\\
  \label{conv-b}
  \mathcal{C}_h^{dg}(w,u,v)=&
  -\frac12\sum_{K\in\Oh}
  \left(\int_{K}(wu)\cdot\nabla v \,\mathrm{dx}
  -\int_{\partial K}\avg{w}u^-\cdot \bld n\, v\,\mathrm{ds}\right),\\
  \label{visc-b}
  \mathcal{B}_h^{hdg}((u,\widehat{u}),
  (v,\widehat{v}))=&\;
  \sum_{K\in\Oh}\Big(\int_{K}\nabla u\cdot\nabla v\,\mathrm{dx}
    -\int_{\partial K}\nabla u\cdot\bld n(v-\widehat{v})\mathrm{ds}\\
  &\;
    -\int_{\partial K}\nabla v\cdot\bld n(u-\widehat{u})\mathrm{ds}
    +
    \int_{\partial K}\frac{4k^2}{h}(u-\widehat{u})(v-\widehat{v})\mathrm{ds}
                 \Big),\nonumber
\end{alignat}
\end{subequations}
where $\avg{w}$ in \eqref{conv-b}  is the standard average operator on element boundaries, and $u^-$ is the {\it upwinding} numerical flux,
with $u^-|_F = (u|_{K^-})|_F$ for any facet $F$ shared by two
elements $K^\pm$, and $K^-$ is the element such that $(\avg{w}\cdot\bld
n_{K^-})|_F\ge 0$.

\subsection{The semidiscrete HDG/DG scheme: Navier-Stokes equations}
The divergence-free HDG/DG scheme in \cite{Lehrenfeld10,LehrenfeldSchoberl16}
is used for the Navier-Stokes equations \eqref{eq:ns}. The semidiscrete scheme reads as follows:
Given initial data $\bld u_h(0)\in \Vh{k}$, for all 
$t\in (0,T]$, find $(\bld u_h,\widehat{\bld u}_h,p_h)=
(\bld u_h(t), \widehat{\bld u}_h(t),p_h)\in \Vh{k}
\times \Vhat{k}\times \Wh{k-1}$ such
that 
\begin{align}
  \label{semi-n}
  \bld{\mathcal{M}}_h(\frac{\partial \bld u_h}{\partial t}, \bld v)
  +\bld{\mathcal{C}}_h^{dg}(\bld u_h,\bld u_h,\bld v)
  +\nu\bld{\mathcal{B}}_h^{hdg}((\bld u_h, \widehat{\bld u}_h),
  (\bld v, \widehat{\bld v}))-\bld{\mathcal{D}}_h(\bld u_h, q)
-\bld{\mathcal{D}}_h(\bld v, p_h)=0,
\end{align}
for all $(\bld v, \widehat{\bld v}, q)\in\Vh{k}\times \Vhat{k}\times
\Wh{k-1}$, where the operators are given as follows: 
\begin{subequations}
  \label{ops-n}  
  \begin{alignat}{2}
  \label{mass-n}
  \bld{\mathcal{M}}_h(\bld u,\bld v)=&\;
  \sum_{K\in\Oh}\int_{K}\bld u\cdot\bld v\,\mathrm{dx},\\
  \label{conv-n}
  \bld {\mathcal{C}}_h^{dg}(\bld w,\bld u,\bld v)=&
  -\sum_{K\in\Oh}
  \left(\int_{K}(\bld w\otimes \bld u):\nabla \bld v \,\mathrm{dx}
  -\int_{\partial K}(\bld w\cdot\bld n)(\bld u^-\cdot \bld v)\,\mathrm{ds}\right),\\
  \label{incomp-n}
  \bld{\mathcal{D}}_h(\bld u,q)=&\;
  \sum_{K\in\Oh}\int_{K}(\nabla\cdot\bld u) q\,\mathrm{dx},\\
  \label{visc-n}
  \bld {\mathcal{B}}_h^{hdg}((\bld u,\widehat{\bld u}),
  (\bld v,\widehat{\bld v}))=&\;
  \sum_{K\in\Oh}\Big(\int_{K}\nabla \bld u:\nabla \bld v\,\mathrm{dx}
    -\int_{\partial K}(\nabla \bld u\bld n)\cdot\mathsf{tang}(\bld v-\widehat{\bld v})\mathrm{ds}\\
                             &\hspace{-6ex}
  -\int_{\partial K}(\nabla \bld v\bld n)\cdot\mathsf{tang}(\bld u-
  \widehat{\bld u})\mathrm{ds}
    +
    \int_{\partial K}\frac{4k^2}{h}
    \mathsf{tang}(\bld u-\widehat{\bld u})\cdot
    \mathsf{tang}(\bld v-\widehat{\bld v})\mathrm{ds}
                 \Big),\nonumber
\end{alignat}
\end{subequations}
where $\bld u^-$ in \eqref{conv-n} is the {\it upwinding} numerical flux,
with $\bld u^-|_F = (\bld u|_{K^-})|_F$ for any facet $F$ shared by two
elements $K^\pm$, and $K^-$ is the element such that $(\bld w\cdot\bld
n_{K^-})|_F\ge 0$, and $\mathsf{tang}(\bld v)|_F:=\bld v-(\bld v\cdot\bld
n)\bld n$ is the tangential component of the vector $\bld v$.
Notice that the convective operator \eqref{conv-n} introduces numerical
dissipation  along element boundaries:
\begin{align}
  \label{conv-stab}
  \bld {\mathcal{C}}_h^{dg}(\bld u_h,\bld u_h,\bld u_h)=&
  \frac12\sum_{K\in\Oh}\int_{\partial K}|\bld u_h\cdot\bld n|
  (\jmp{\bld  u_h})^2\,\mathrm{ds}\ge0, \quad \forall \bld u_h\in \Vh{k}\cap
  H(\mathrm{div}^0,\Omega),
\end{align}
where 
\[
  H(\mathrm{div}^0,\Omega):=\{\bld v\in H(\mathrm{div},\Omega):\;\;\nabla
  \cdot\bld v=0 \},
\] 
which is beneficial in the convection-dominated regime.

We remark that the scheme \eqref{semi-n} produces an exactly divergence-free
velocity approximation, i.e. 
$\bld u_h\in \Vh{k,0}:=\Vh{k}\cap H(\mathrm{div}^0,\Omega)$, which is a desired property 
for the POD model we consider in the next section. 
In particular, the velocity field 
$(\bld u_h, \widehat{\bld u}_h)\in \Vh{k,0}\times \Vhat{k}$
can be directly computed without pressure approximation by solving the
following equations:
\begin{align}
  \label{semi-n2}
  \bld{\mathcal{M}}_h(\frac{\partial \bld u_h}{\partial t}, \bld v)
  +\bld{\mathcal{C}}_h^{dg}(\bld u_h,\bld u_h,\bld v)
  +\nu\bld{\mathcal{B}}_h^{hdg}((\bld u_h, \widehat{\bld u}_h),
  (\bld v, \widehat{\bld v}))
  =0,\quad\forall
(\bld v, \widehat{\bld v})\in \Vh{k,0}\times \Vhat{k}.
\end{align}

\subsection{The fully discrete HDG schemes}
For the time discretization, we use
the second-order Crank-Nicolson-Adams-Bashforth (CNAB) method \cite{AscherRuuthWetton93}, which treats
the nonlinear convective term explicitly, and other terms implicitly.
For simplicity, a uniform time partition is applied. Let $0=t_0< t_1<\cdots<t_M=T$ be the partition of the interval $[0,T]$ and the  
time step $\Delta t = \frac{T}{M}$. 

The fully discrete scheme for the Burgers' equation \eqref{eq:burgers} is given as follows:
Given initial data $(u_h^0, \widehat{u}_h^0)\in \Wh{k}\times 
\What{0}$, for each integer $n=1,\cdots, M$, find
$(u_h^n,
\widehat{u}_h^{n})\in \Wh{k}\times \What{0}$ such that 
\begin{align}
  \label{full-b}
  \mathcal{M}_h(\frac{u_h^{n}-u_h^{n-1}}{\Delta t}, v)
  +\mathcal{C}_h^{dg}(\widetilde{u}_h^{n-1/2},
  \widetilde{u}_h^{n-1/2},v)
  +\nu\mathcal{B}_h^{hdg}((u_h^{n-1/2}, 
  \widehat{u}_h^{n-1/2}),
  (v, \widehat{v})) = 0, 
\end{align}
for all $(v,\widehat{v})\in \Wh{k}\times \What{0}$,
  where 
  \[
    u_h^{n-1/2}:=\frac12(u_h^{n}+u_h^{n-1}),\quad 
    \widehat{u}_h^{n-1/2}:=
    \frac12(\widehat{u}_h^{n}+
    \widehat{u}_h^{n-1}),\quad 
    \widetilde{u}_h^{n-1/2}:=\frac32u_h^{n-1}-\frac12u_h^{n-2}.
  \]
  Here in the first step ($n=1$)  we simply take $\widetilde{u}_h^{1/2}=u_h^{0}$. 

Similarly,  
the fully discrete scheme for the Navier-Stokes equations
\eqref{eq:ns}  is given as follows:
  Given initial data $(\bld u_h^0,\widehat{\bld u}_h^0)\in \Vh{k}
  \times \Vhat{k}$, for each integer $n=1,\cdots, M$, find
$(\bld u_h^n,
\widehat{\bld u}_h^{n},p_h^{n-1/2})\in \Vh{k}\times \Vhat{k}\times
\Wh{k-1}$ such that 
\begin{align}
  \label{full-n}
  &
  \bld{\mathcal{M}}_h(\frac{\bld u_h^{n}-\bld u_h^{n-1}}{\Delta t}, 
  \bld v)
  +\bld{\mathcal{C}}_h^{dg}(\widetilde{\bld u}_h^{n-1/2},
  \widetilde{\bld u}_h^{n-1/2},\bld v)
  +\nu\bld{\mathcal{B}}_h^{hdg}((\bld u_h^{n-1/2}, 
  \widehat{\bld u}_h^{n-1/2}),
  (\bld v, \widehat{\bld v})) \\
  &\hspace{44ex} -\bld{\mathcal{D}}_h(\bld v, p_h^{n-1/2})
 -\bld{\mathcal{D}}_h(\bld u_h^{n-1/2},q)  = 0,\nonumber 
\end{align}
for all $(\bld v,\widehat{\bld v},q)\in \Vh{k}\times \Vhat{k}\times
\Wh{k-1}$.

Efficient implementation of the HDG linear system \eqref{full-b} and
\eqref{full-n} via static condensation were discussed, for example, in \cite{Cockburn16,Lehrenfeld10}.

  \section{The POD model}
  \label{sec:pod}
In this section, we present the POD model based on the FOM IMEX HDG/DG
schemes presented in
Section \ref{sec:fom} using the method of snapshots \cite{Sirovich87}. 
We focus on the discussion for the 
Navier-Stokes equations as the results for Burgers' equation are identical.
 Since the generated POD basis functions are {\it global}, 
 we do not see any 
 advantage of formulating a POD-HDG 
ROM
 constructed
 using both variables $\bld u_h$ and $\widehat{\bld u}_h$.
Hence,  we only use the field variable $\bld u_h$ to construct the POD
model, and the resulting ROM is a DG scheme.

\subsection{Computing POD basis functions}
The method of snapshots is used to construct the POD bases.
To this end, 
let $\{\bld u_h^n\}_{n=0}^{S-1}$ be snapshots obtained from a full order model
simulation \eqref{full-n}.
The POD bases are obtained by the following steps:
\begin{itemize}
  \item [(i)]
  Decompose the data $\bld u_h^{n}$ into the mean part ($\bar{\bld u}_h$) and the fluctuating part
($\check{\bld u}_h^n$):
\[
  \bld u_h^n = \bar{\bld u}_h + \check{\bld u}_h^n,
  \quad\bar{\bld u}_h=\frac1S\sum_{n=0}^{S-1}\bld u_h^n.
\]
\item [(ii)]
  Build the (symmetric positive definite) correlation matrix $C\in
\mathbb{R}^{S\times S}$ with
$
  C_{ij} = \bld{\mathcal{M}}_h(\check{\bld u}_h^i,\check{\bld u}_h^j).
$
\item [(iii)] Solve the eigenvalue problem:
  \[
  CW = W\Lambda,
  \] 
  where $\Lambda=\mathrm{diag}[\lambda_1,\cdots,\lambda_S]$, 
  $W=[w^1,\cdots,w^S]$, $\lambda_i$ is the $i$th eigenvalue and $w^i$ is
  the corresponding normalized $i$th eigenvector.
\item [(iv)] Given an integer $r\ll S$, 
  return the first $r$ POD basis functions $\{\phi_j\}_{j=1}^r$, where
  \[
    \bld \phi_j
    = \frac{1}{\sqrt{\lambda_j}}\sum_{n=0}^{r-1}w_n^j\check{\bld u}_h^n,
    \quad j=1,\cdots, r.
  \] 
\end{itemize}
Denote the space $S_h^r=\mathrm{span}\{\bld \phi_1,\cdots\bld\phi_r\}$.    
Since ${\bld \phi_j}$ are orthonormal, the mass matrix associated with the space
$S_h^r$ is the identity matrix.

\subsection{The plain POD-DG scheme}
To construct the POD-DG scheme, we first replace the HDG viscous operator
\eqref{visc-n}
by a DG operator:
\begin{align}
  \label{visc-dg}
  \bld{\mathcal{B}}_h^{dg}(\bld u, \bld v)=&\;
  \sum_{K\in\Oh}\Big(\int_{K}\nabla \bld u:\nabla \bld v\,\mathrm{dx}
    -\int_{\partial K}(\avg{\nabla \bld u}\bld n)\cdot\mathsf{tang}(\bld v)\mathrm{ds}\\
                             &\hspace{-6ex}
                             -\int_{\partial K}(\avg{\nabla \bld v}\bld
                             n)\cdot\mathsf{tang}(\bld u)\mathrm{ds}
    +
    \int_{\partial K}\frac{4k^2}{h}
    \mathsf{tang}(\jmp{\bld u})\cdot
    \mathsf{tang}(\jmp{\bld v})\mathrm{ds}
                 \Big),\nonumber
\end{align}
where, on each internal facet $F\in\Eh^i$, 
$\jmp{\bld v}|_F=(\bld v^+-\bld v^-)$ is the standard jump operator, and 
$\jmp{\bld u}|_{\partial\Omega}=0$.
Next, we notice that the upwinding convection operator \eqref{conv-n} is
linear in the first and third arguments, but nonlinear in the second
argument, due to the upwinding numerical flux $\bld u^-$.
This nonlinearity is quite troublesome for ROM in the sense that it can not be computed
using an offline procedure. 
We mention that it is precisely this nonlinear term that provides the
upwinding mechanism for the DG operator \eqref{conv-n}, 
which produces extra numerical dissipation to stabilize the scheme
\eqref{full-n} in the under-resolved convection-dominated
regime.
To seek for an efficient implementation, we replace 
the ({\it nonlinear}) 
upwinding flux by  the ({\it linear}) central flux:
\begin{align}
  \label{conv-c}
  \bld {\widetilde{\mathcal{C}}}_h^{dg}(\bld w,\bld u,\bld v)=&
  -\frac12\sum_{K\in\Oh}
  \left(\int_{K}(\bld w\otimes \bld u):\nabla \bld v \,\mathrm{dx}
  -\int_{\partial K}(\bld w\cdot\bld n)(\avg{\bld u}\cdot \bld
v)\,\mathrm{ds}\right).
\end{align}
This creates a {\it trilinear} operator that satisfies
the following energy conservation property:
\begin{align}
  \label{conv-conv}
  \bld {\widetilde{\mathcal{C}}}_h^{dg}(\bld u_h,\bld u_h,\bld u_h)=&
  0, \quad \forall \bld u_h\in \Vh{k}\cap
  H(\mathrm{div}^0,\Omega).
\end{align}
Finally, the pressure field can be directly eliminated from the 
POD scheme because all POD basis functions are globally divergence-free, inherited from the snapshots.
The semidiscrete plain POD-DG scheme reads as follows:
Given initial data $\bld u_h(0)=\bar{\bld u}_h+\check{\bld u}_h^0$ with
$\check{\bld u}_h^0\in S_h^r$, for all 
$t\in (0,T]$, find $\bld u_h=\bar{\bld u}_h+\check{\bld u}_h(t)$ with
$\check{\bld u}_h(t)\in S_h^r$ such
that 
\begin{align}
  \label{semi-pod}
  \bld{\mathcal{M}}_h(\frac{\partial \bld u_h}{\partial t}, \bld v)
  +\bld{\widetilde{\mathcal{C}}}_h^{dg}(\bld u_h,\bld u_h,\bld v)
  +\nu\bld{\mathcal{B}}_h^{dg}(\bld u_h,\bld v) = 0, \quad\forall \bld v\in S_h^r.
\end{align}

We again use the CNAB time discretization, and the fully discrete plain
POD-DG scheme reads as follows:
Given initial data $\bld u_h^0=\bar{\bld u}_h+\check{\bld u}_h^0$ with
$\check{\bld u}_h^0\in S_h^r$, for each integer $n=1,\cdots, M$, find
$\bld u_h^n=\bar{\bld u}_h+\check{\bld u}_h^n$ with $\check{\bld u}_h^n\in S_h^r$ such that 
\begin{align}
  \label{full-pod}
  \bld{\mathcal{M}}_h(\frac{\bld u_h^{n}-\bld u_h^{n-1}}{\Delta t}, 
  \bld v)
  +\bld{\widetilde{\mathcal{C}}}_h^{dg}(\widetilde{\bld u}_h^{n-1/2},
  \widetilde{\bld u}_h^{n-1/2},\bld v)
  +\nu\bld{\mathcal{B}}_h^{dg}(\bld u_h^{n-1/2}, \bld v) =0,\quad
\forall \bld v\in S_h^r.
\end{align}

\begin{remark}[Offline-online decomposition]
  The POD-DG schemes \eqref{semi-pod} and  \eqref{full-pod} can be efficiently implemented via
  a standard offline-online decomposition.
  Thus, we introduce the offline-precomputable vectors 
  $\mathsf{C0}, \mathsf{B0}\in
  \mathbb{R}^r$, matrices $\mathsf{C1}, \mathsf{B}\in
  \mathbb{R}^{r\times r}$, and third
  order tensor $\mathsf{C}\in \mathbb{R}^{r\times r\times r}$:
    \begin{align*}
      \mathsf{C0}_j=&\; \bld{\widetilde{\mathcal{C}}}_h^{dg}(\bar{\bld
      u}_h,\bar{\bld u}_h, \bld \phi_j),\\
      \mathsf{B0}_j=&\; \bld{\mathcal{B}}_h^{dg}(\bar{\bld u}_h, \bld \phi_j),\\ 
      \mathsf{C1}_{i,j}=&\;
      \bld{ \widetilde{\mathcal{C}}}_h^{dg}(\bar{\bld u}_h,\bld\phi_i,\bld \phi_j)
      +
      \bld{\widetilde{\mathcal{C}}}_h^{dg}(\bld\phi_i, \bar{\bld u}_h, \bld\phi_j),\\
      {\mathsf{B}}_{i,j}=&\;\bld{ \mathcal{B}}_h^{dg}(\bld\phi_i,\bld \phi_j),\\ 
      { \mathsf{C}}_{i,j,k}=&\;
      \bld{\widetilde{\mathcal{C}}}_h^{dg}(\bld \phi_i,\bld\phi_j, \bld\phi_k).
    \end{align*}
    Denote $\check{ \bld u}_h=\sum_{j=1}^r a_j(t)\bld\phi_j$, then the semi-discrete
  scheme \eqref{semi-pod} is given in the following form:
  \begin{align}
    \label{semi-pod-online}
    \frac{\partial a_j}{\partial t}
    +\mathsf{C0}_j+\mathsf{C1}_{ij}a_i+
    \mathsf{C}_{ikj}a_ia_k
    +\nu\mathsf{B0}_j+\nu\mathsf{B}_{ij}a_i=0
  \end{align}
  Denoting $\bld a^n = [a_1^n,\cdots,a_r^n]\in \mathbb{R}^r$, 
  the fully discrete scheme \eqref{full-pod} 
  is then given in the following form, which can be be computed efficiently online,
  \begin{align}
    \label{full-pod-online}
    (\frac{Id}{\Delta t}
    +\frac12\nu\mathsf{B})\bld a^{n}
    = 
    \bld a^{n-1}    -(\mathsf{C0}+\mathsf{C1}\widetilde{\bld a}^{n-1/2}
    +
  \widetilde{\bld a}^{n-1/2}\mathsf{C}\widetilde{\bld a}^{n-1/2}
    +\nu\mathsf{B0}+\frac12\nu\mathsf{B}\bld a^{n-1}),
  \end{align}
  where $\widetilde{\bld a}^{n-1/2}=\frac32\bld a^{n-1}-\frac12\bld
  a^{n-2}$.
\end{remark}

\subsection{The closure model}
Due to the use of {\it linear} central numerical flux 
for the convection operator, the plain POD-DG scheme \eqref{full-pod} does
not inherit the extra ({\it upwinding}) 
convective stabilization property of the original HDG/DG scheme
\eqref{full-n} that, however, is the key for the stability of the scheme in the
under-resolved convection dominated regime.
Hence, it is natural to introduce a {\it linear} stabilization term 
that mimics such upwinding mechanism in the POD setting.
We further include a standard eddy viscosity closure model originally proposed
in \cite{SanIliescu13}, in order to improve accuracy/stability of
the POD-DG scheme. 
To this end, we denote the following non-negative matrices $\mathsf{CX}, \mathsf{BX}\in \mathbb{R}^{r\times
r}$:
\begin{subequations}
  \label{stab}  
\begin{align}
  \mathsf{CX}_{ik} = &\;
  \frac12\sum_{K\in\Oh}\int_{\partial
  K}\jmp{\bld\phi_i}\cdot\jmp{\bld\phi_k}\,\mathrm{ds},\\
  \mathsf{BX}_{ik} = &\; \bld{\mathcal{B}}_h^{dg}(\bld\phi_i, (k/r)^2
  \bld \phi_k),
\end{align}
\end{subequations}
and define the POD-DG closure model as follows:
  \begin{align}
  \label{closure}
    (\frac{Id}{\Delta t}
    +\frac12\widetilde{\mathsf{B}})\bld a^{n}
    = 
    \bld{a}^{n-1}-(\mathsf{C0}+\mathsf{C1}\widetilde{\bld a}^{n-1/2}
    +
  \widetilde{\bld a}^{n-1/2}\mathsf{C}\widetilde{\bld a}^{n-1/2}
  +\nu\mathsf{B0}+\frac12\widetilde{\mathsf{B}}\bld a^{n-1}),
  \end{align}
  where $\widetilde{\mathsf{B}}=\nu\mathsf{B}+c_1\mathsf{CX}
  +c_2\mathsf{BX}$, with $c_1,c_2\ge 0$ being two tunable constants that
  are
  problem dependent.
  Here the matrix $\mathsf{BX}$ corresponds to an eddy viscosity model with 
a quadratic viscosity kernel \cite{SanIliescu13},  and the matrix 
$\mathsf{CX}$ can be interpreted as an upwinding 
stabilization term (compared with DG upwinding in \eqref{conv-stab}).
We call the term with $\mathsf{BX}$ a diffusive stabilization, and the term
with $\mathsf{CX}$ a convective stabilization.
We remark that if we take $c_1=\max{|\bld u_h|}$, then 
the convective stabilization term scales similarly as the full order model
case. 
However, our numerical results in the next section indicates that
taking $c_1=\max{|\bld u_h|}$ is too small to make such convective
stabilization term effective in the POD setting. Actually, in  a case for
the Burgers' equation, 
we need to take $c_1=2\times 10^{8}$ (see Example 1 in Section \ref{sec:num} below) to see the positive impact of
this stabilization term. 
This observation also partially justify our choice of 
{\it linear} central numerical flux in the convection operator 
\eqref{conv-c}
over the {\it nonlinear} upwinding
numerical flux for the plain POD-DG scheme \eqref{full-pod}. 
 
Finally,  we remark that the two 
parameters $c_1$ and $c_2$ are tuned purely at the online stage, such tuning cost is negligible comparing
 to the computational cost of the full order model \eqref{full-b}.

\section{Numerical results}
  \label{sec:num}
In this section, we present some numerical examples for the 
POD-DG closure model \eqref{closure}.
The NGSolve software \cite{Schoberl16} is used for the simulations.

\subsection{Example 1. Burgers' equation: discontinuous initial condition}
\begin{figure}[!ht]
\centering
\includegraphics[width=.45\textwidth]{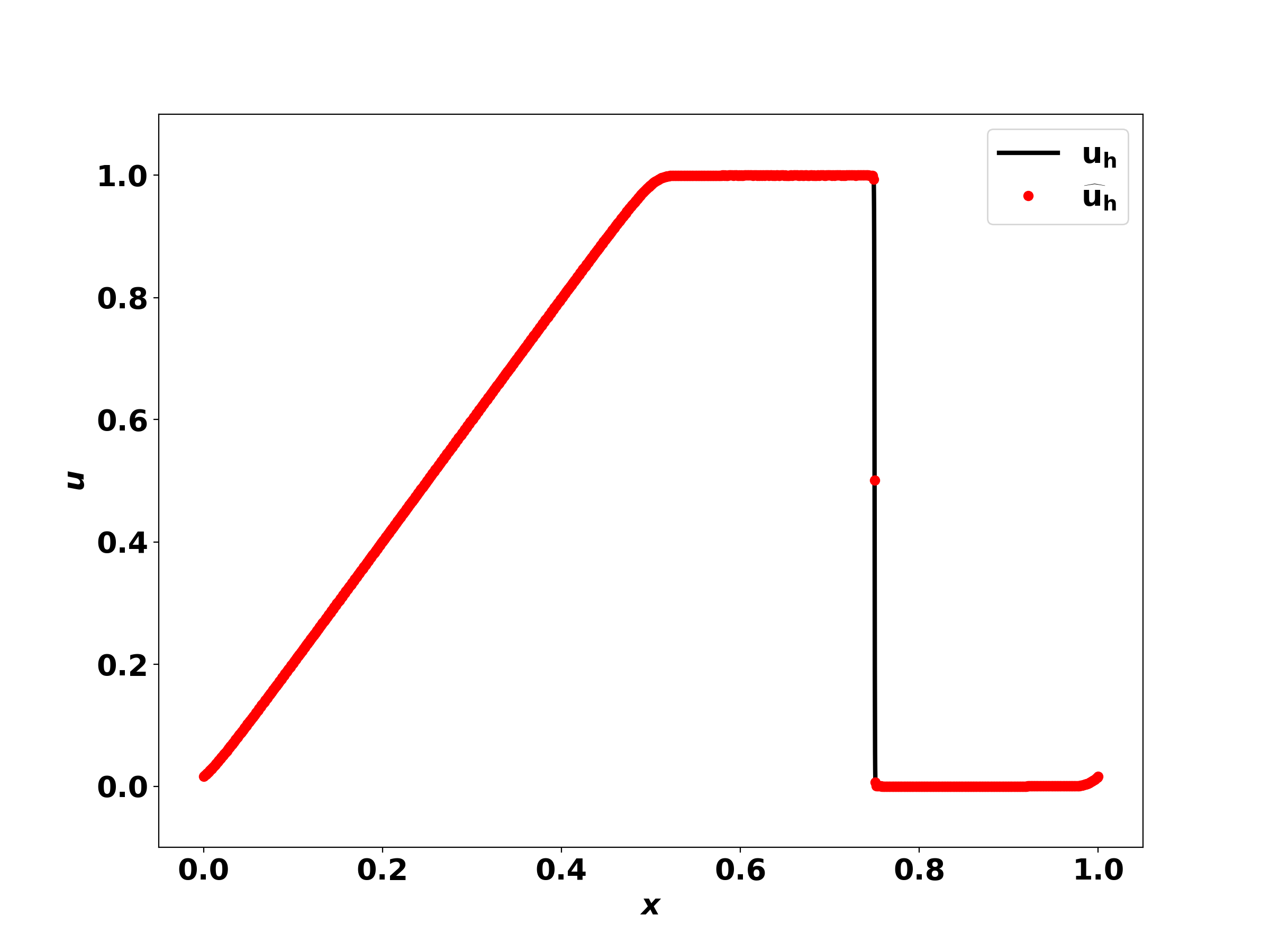}
\includegraphics[width=.45\textwidth]{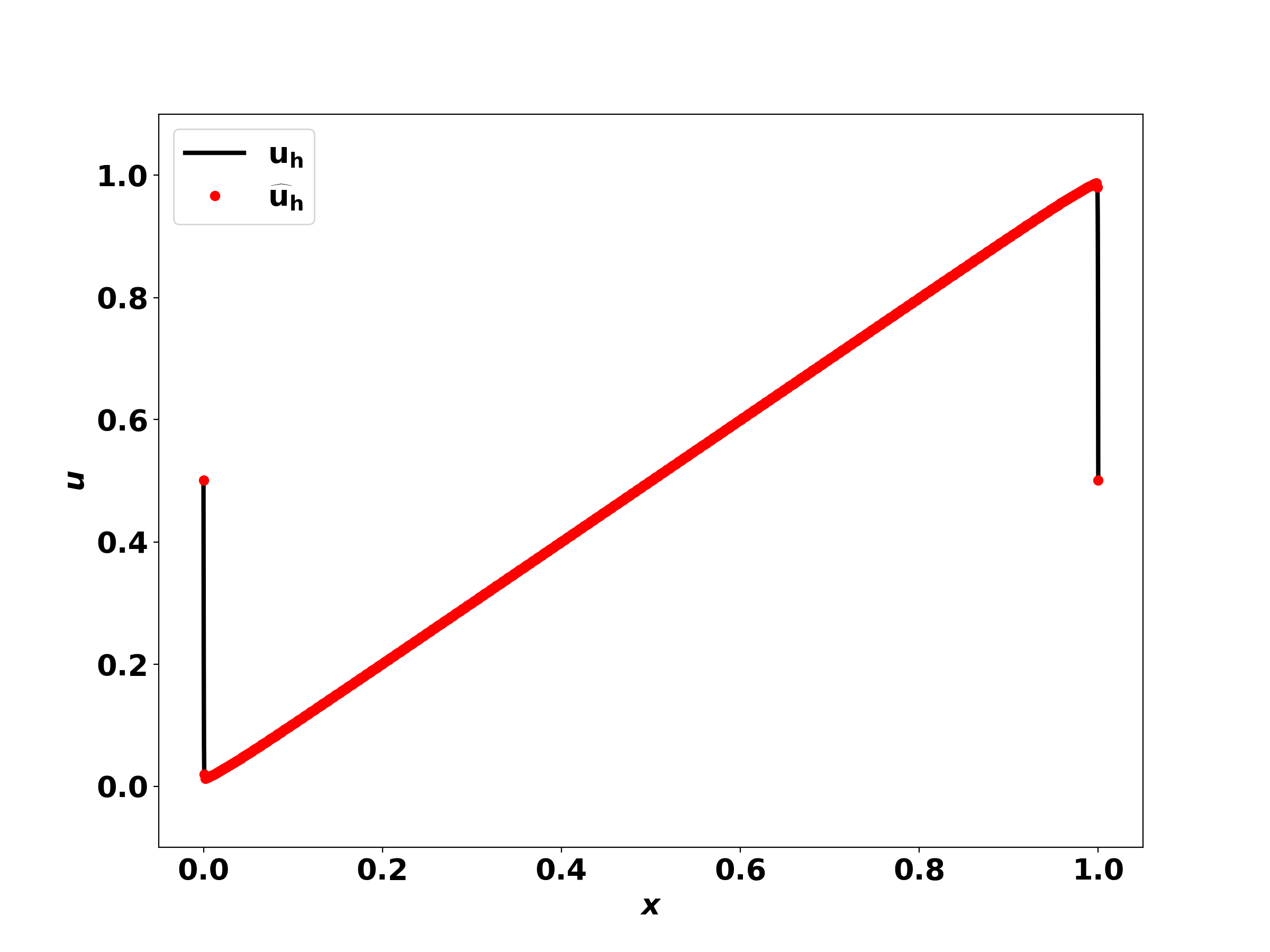}
  \caption{
Example 1.    Numerical solution 
at $t=0.5$ (left) and $t=1$ (right). 
Black line: $u_h$. Red dots: $\widehat{u}_h$.
  }
\label{fig:ex1}
\end{figure}

\begin{figure}[!ht]
\centering
\includegraphics[width=.45\textwidth]{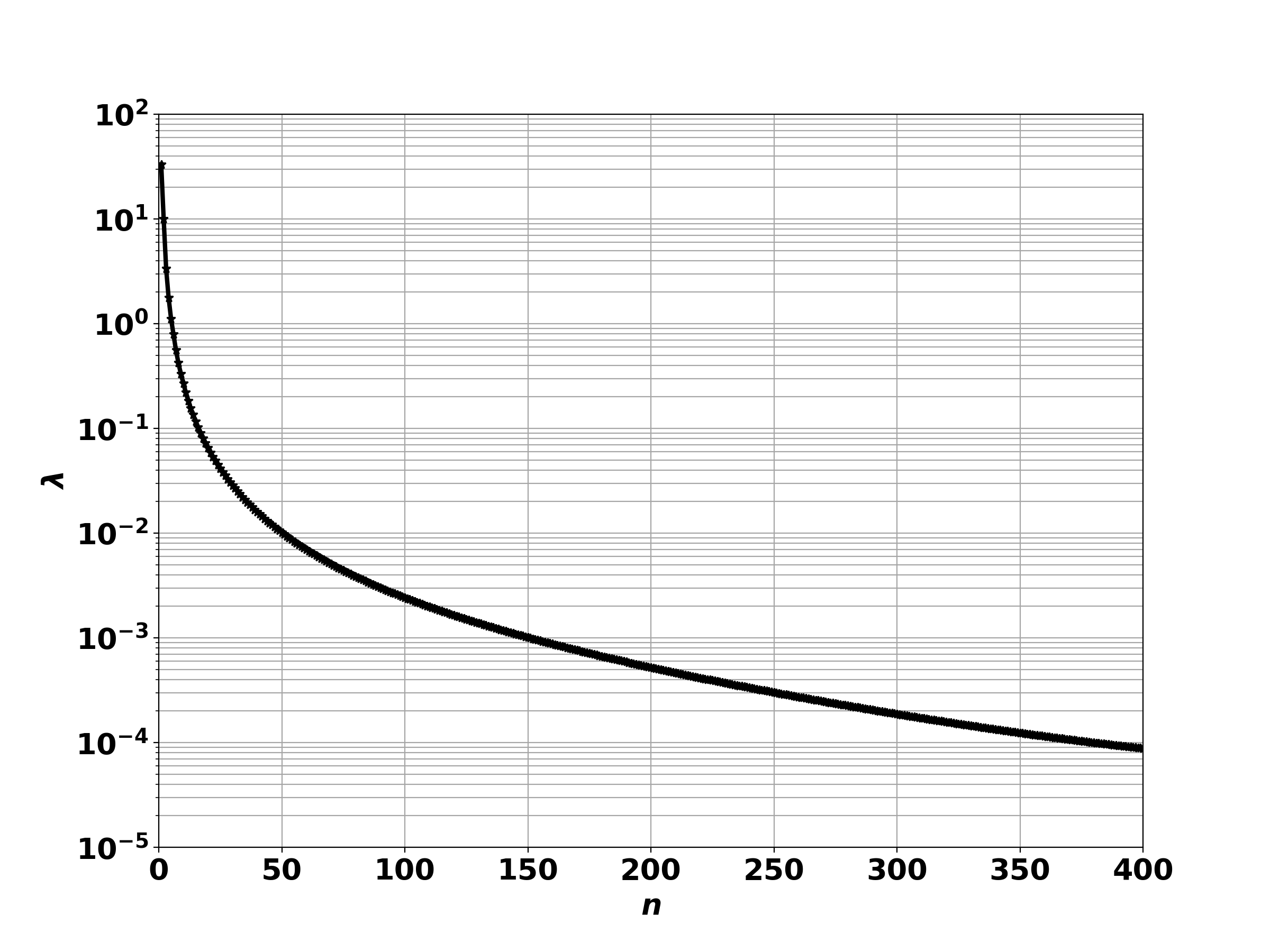}
\includegraphics[width=.45\textwidth]{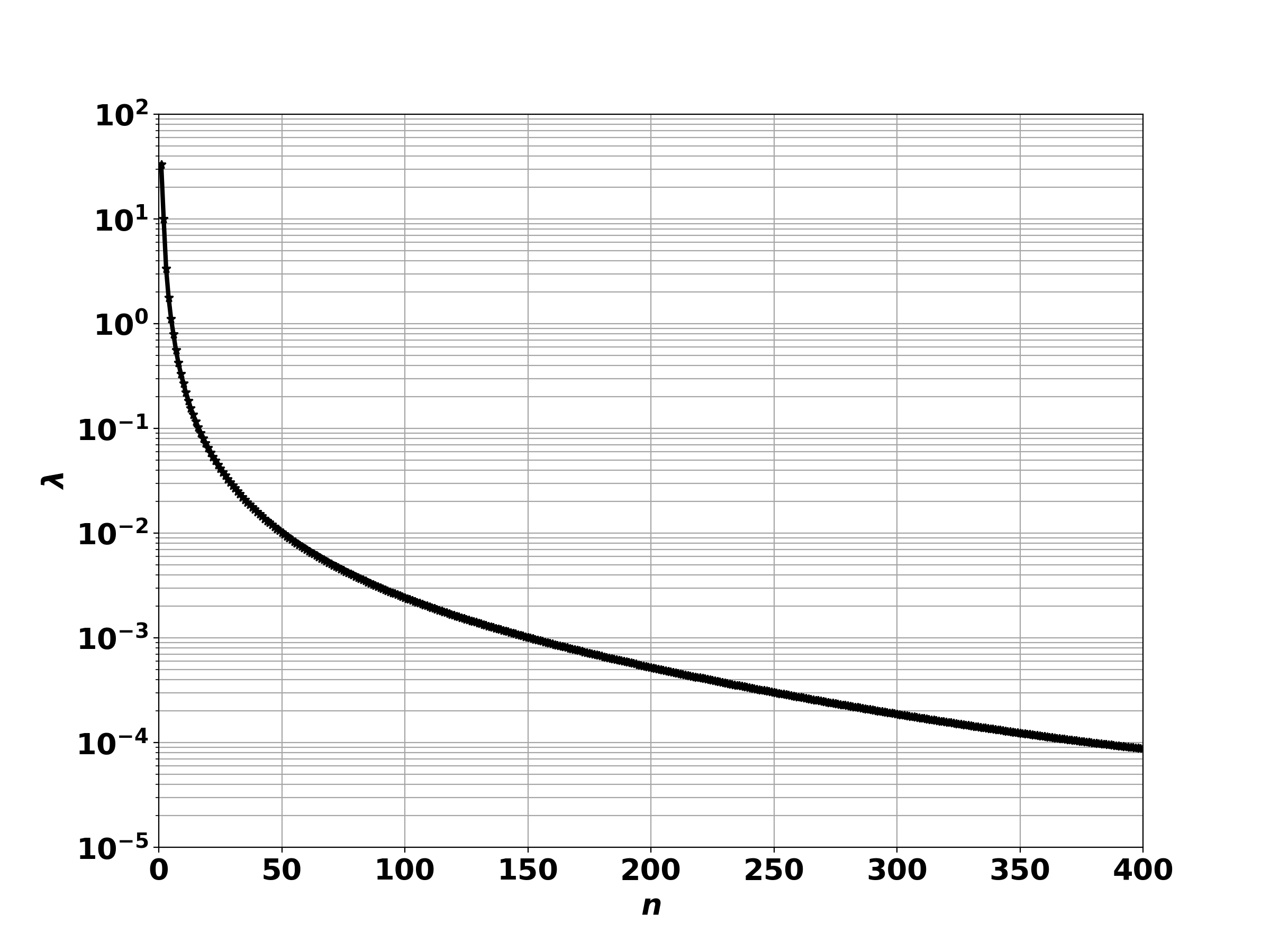}
  \caption{
Example 1. First 400 eigenvalues of the correlation matrix $C$.
    Left: slightly resolved case $k=2$. Right: fully resolved case $k=6$.
  }
\label{fig:ex1-eig}
\end{figure}

\begin{figure}[!ht]
\centering
\includegraphics[width=.45\textwidth]{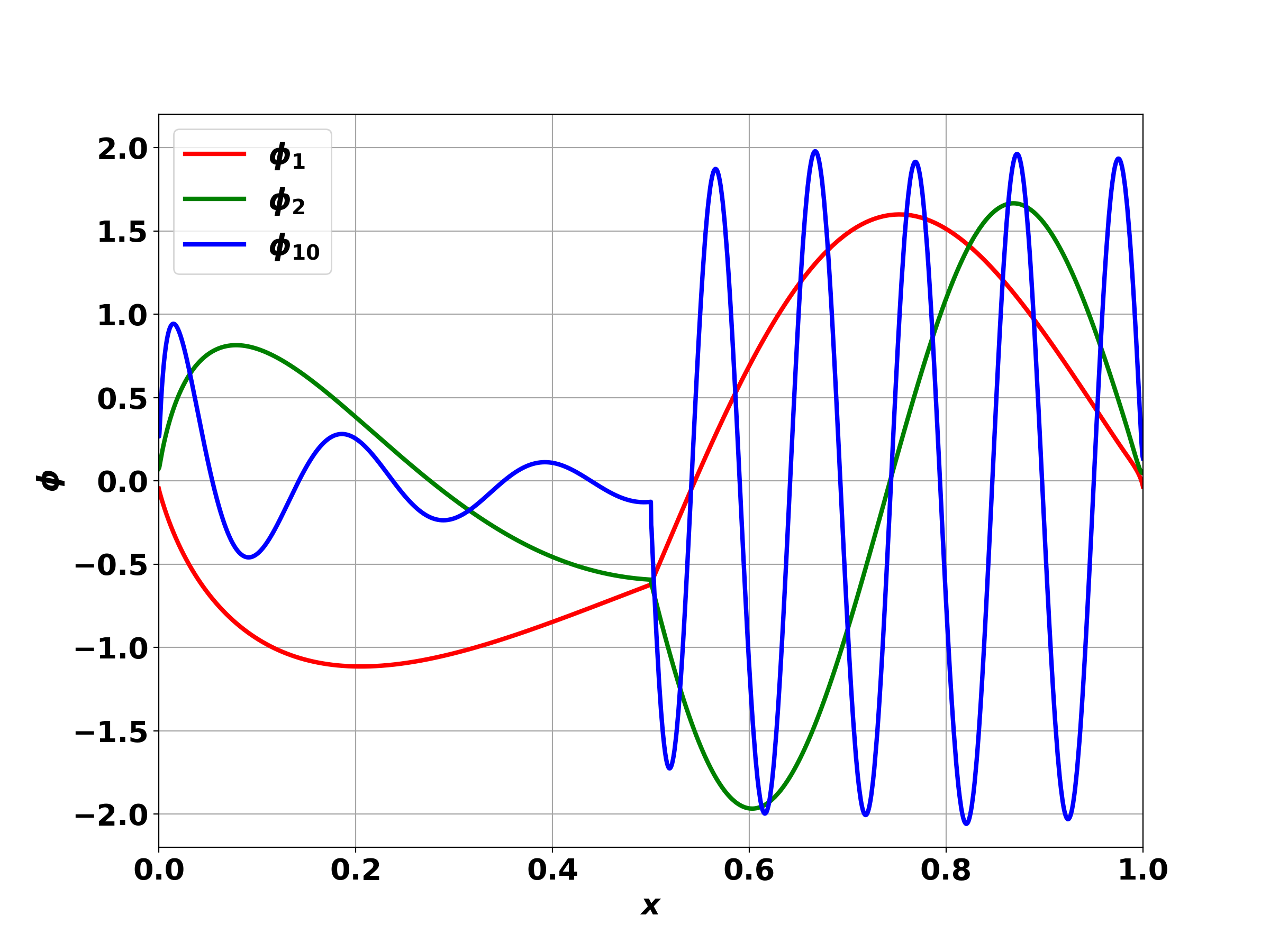}
\includegraphics[width=.45\textwidth]{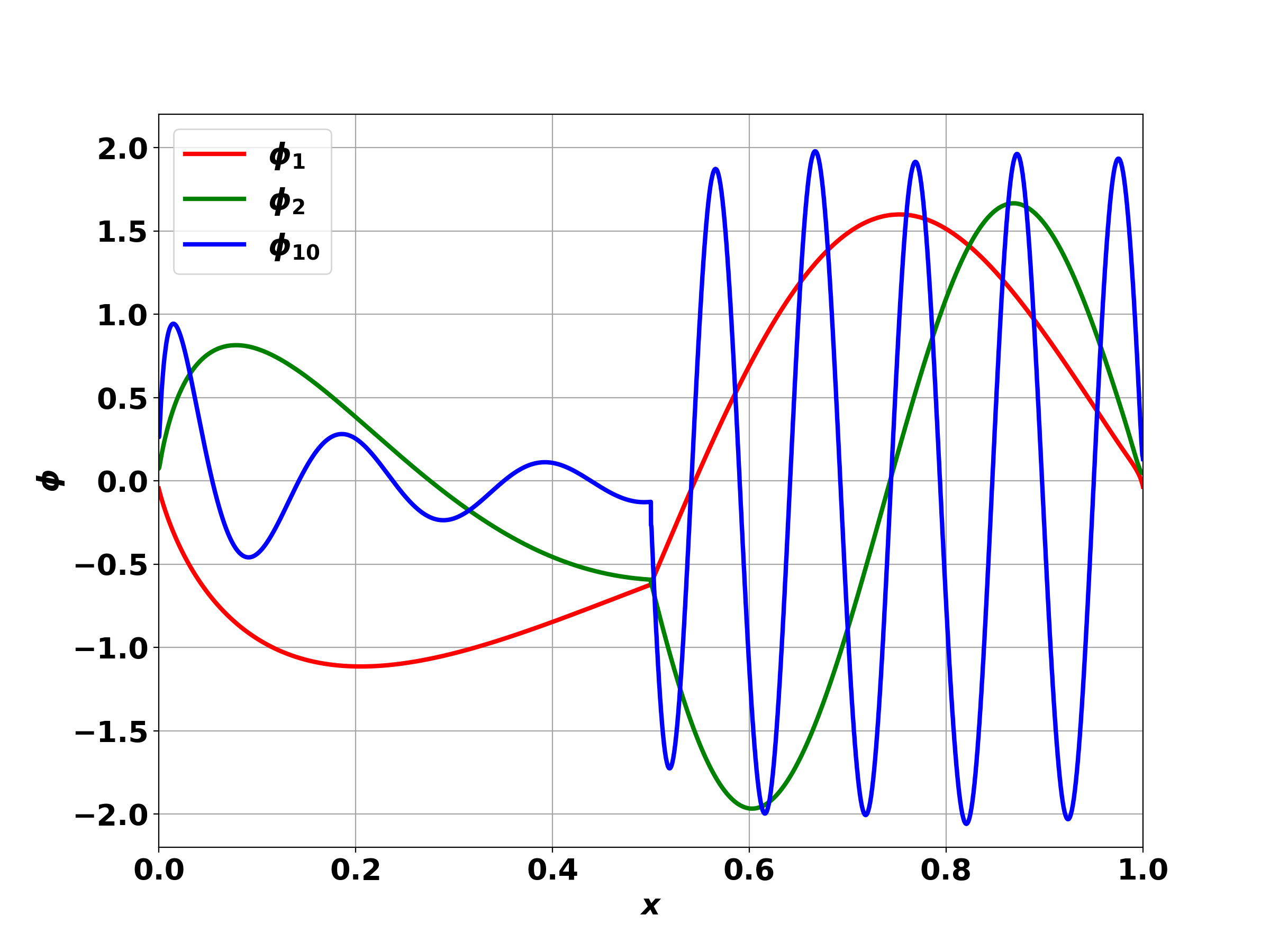}
  \caption{
Example 1.
    Illustrative examples of POD basis functions.
    Left: slightly resolved case $k=2$. Right: fully resolved case $k=6$.
  }
\label{fig:ex1-basis}
\end{figure}

We consider the Burgers' equation \eqref{eq:burgers} with $\nu=10^{-4}$ and the periodic
boundary conditions.
The initial condition is taken to be a step function 
  \[
    u(0) = \left\{
      \begin{tabular}{ll}
        1 & if $x<0.5$,\\[1ex]
        0 & if $x\ge 0.5$.
    \end{tabular}\right.,
  \] 
 and the final time is $T=1$.
Two cases of the full order model \eqref{full-b} are tested here that associate with different discretization parameters, including mesh size $h$ of the uniform mesh, polynomial degree $k$, and the uniform time step size $\Delta t$.  
\begin{itemize}
  \item [(i)] Slightly resolved case:
    $
      h = 10^{-4}, k=2, \Delta t = 0.1h.
    $ 
  \item [(ii)] Fully resolved case:
    $
      h = 10^{-4}, k=6, \Delta t = 0.04h.
    $
\end{itemize}
To build the POD basis, we collect $501$ snapshots in the time interval
$[0,1]$ taken at equidistant time instances.
  The numerical solutions for the fully resolved case ($k=6$) at $t=0.5$
  and $t=1$ are shown in
  Figure~\ref{fig:ex1}. We observe the sharp gradient is resolved within
  2 cells. 
The eigenvalues of the correlation matrix $C$ are shown in
Figure~\ref{fig:ex1-eig}, where we do not observe any significant
difference for both cases.
Figure~\ref{fig:ex1-basis} visualizes three POD basis functions
$\phi_1(x)$, $\phi_2(x)$, and $\phi_{10}(x)$ for both cases. Again, we observe no
significant difference between these two cases.

\begin{figure}[!ht]
\centering
\includegraphics[width=.45\textwidth]{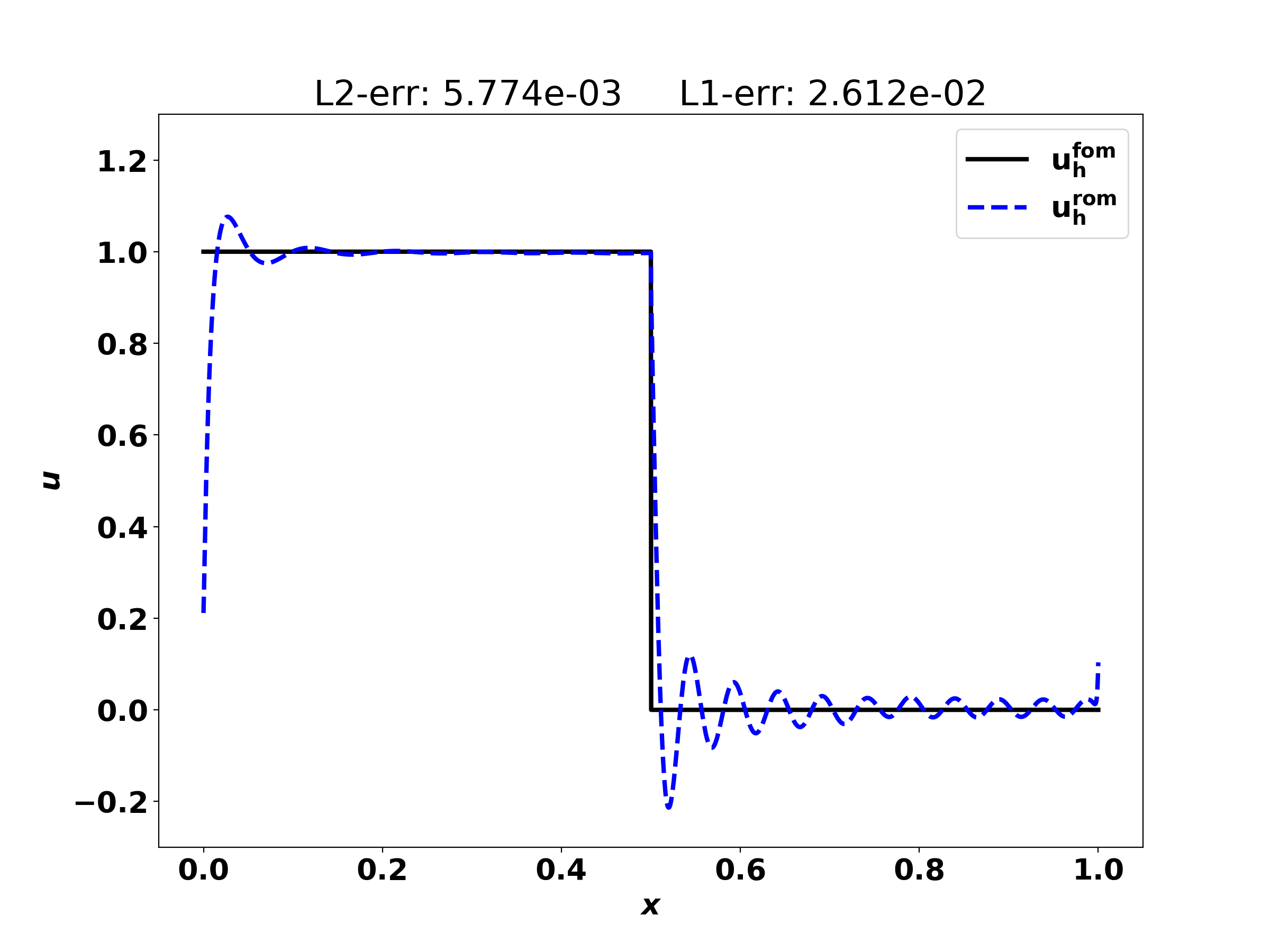}
\includegraphics[width=.45\textwidth]{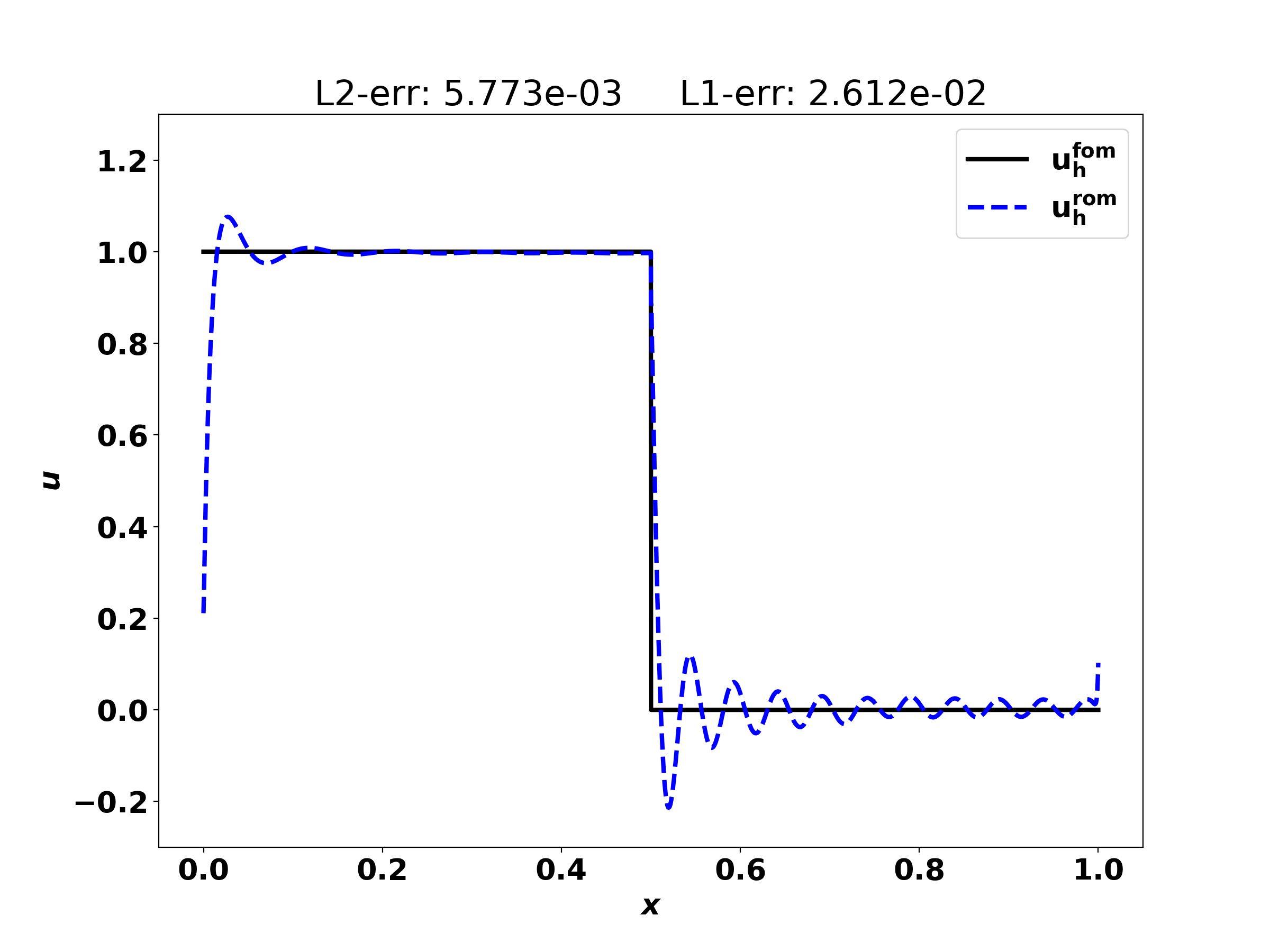}
\includegraphics[width=.45\textwidth]{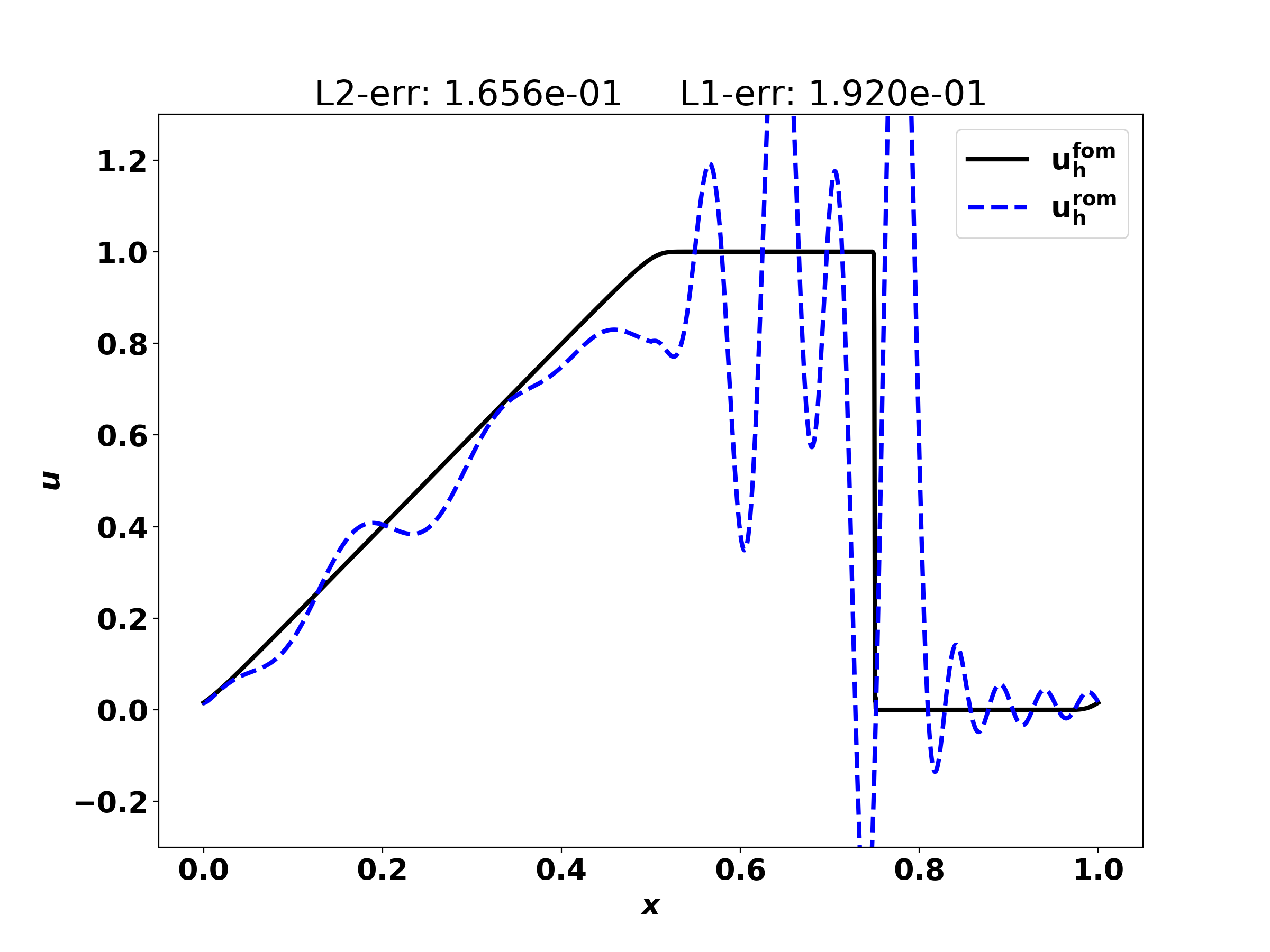}
\includegraphics[width=.45\textwidth]{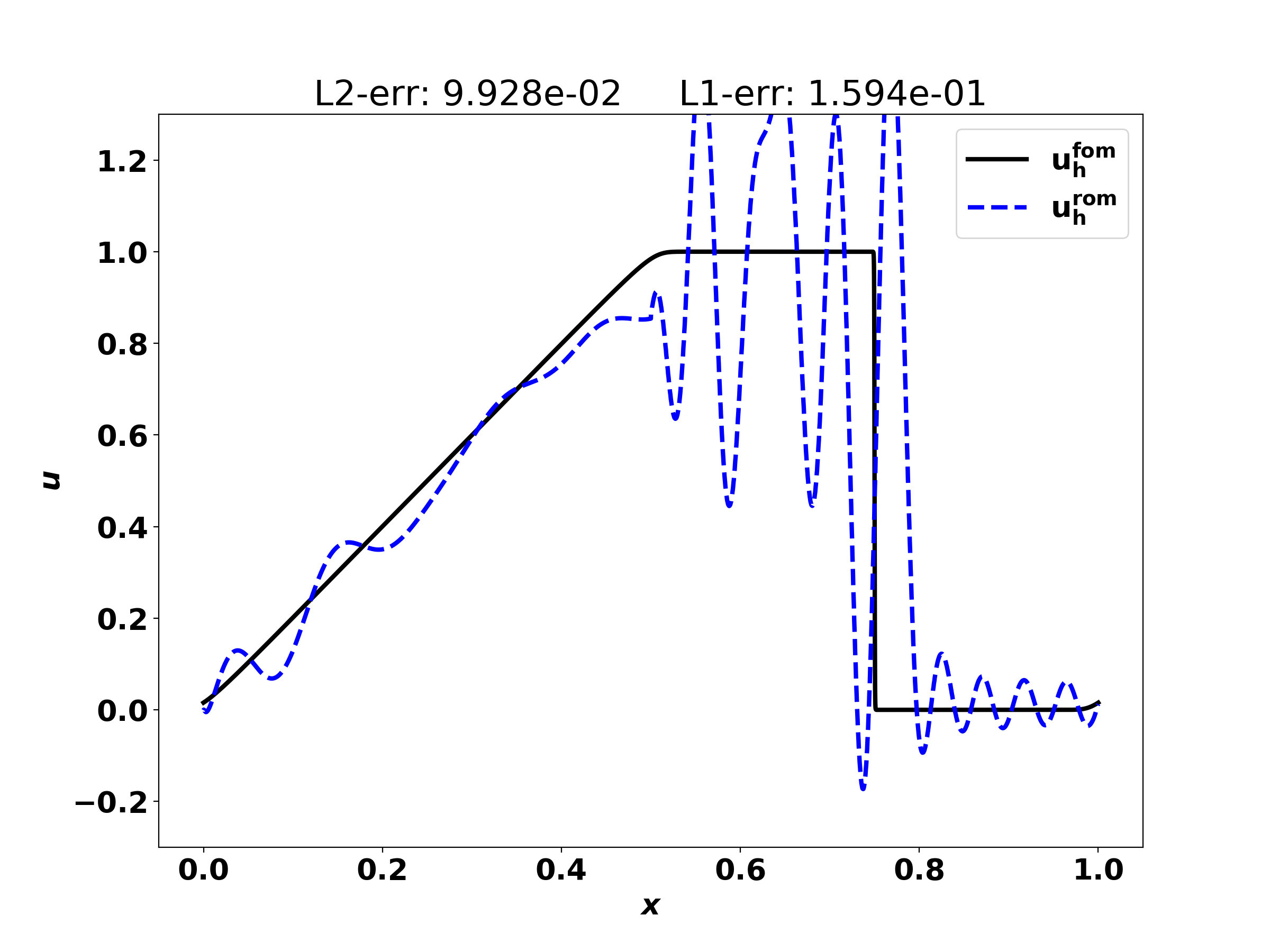}
\includegraphics[width=.45\textwidth]{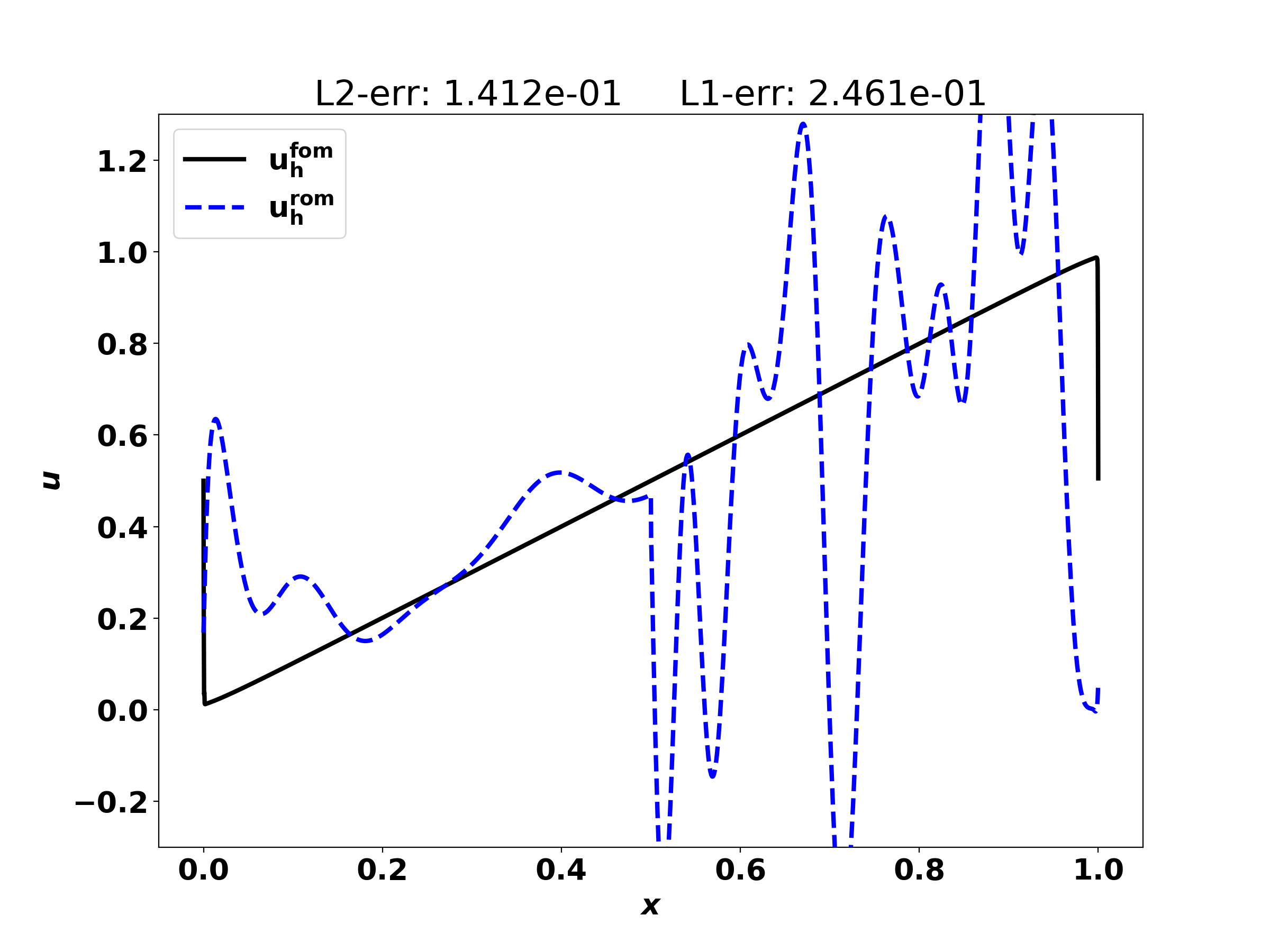}
\includegraphics[width=.45\textwidth]{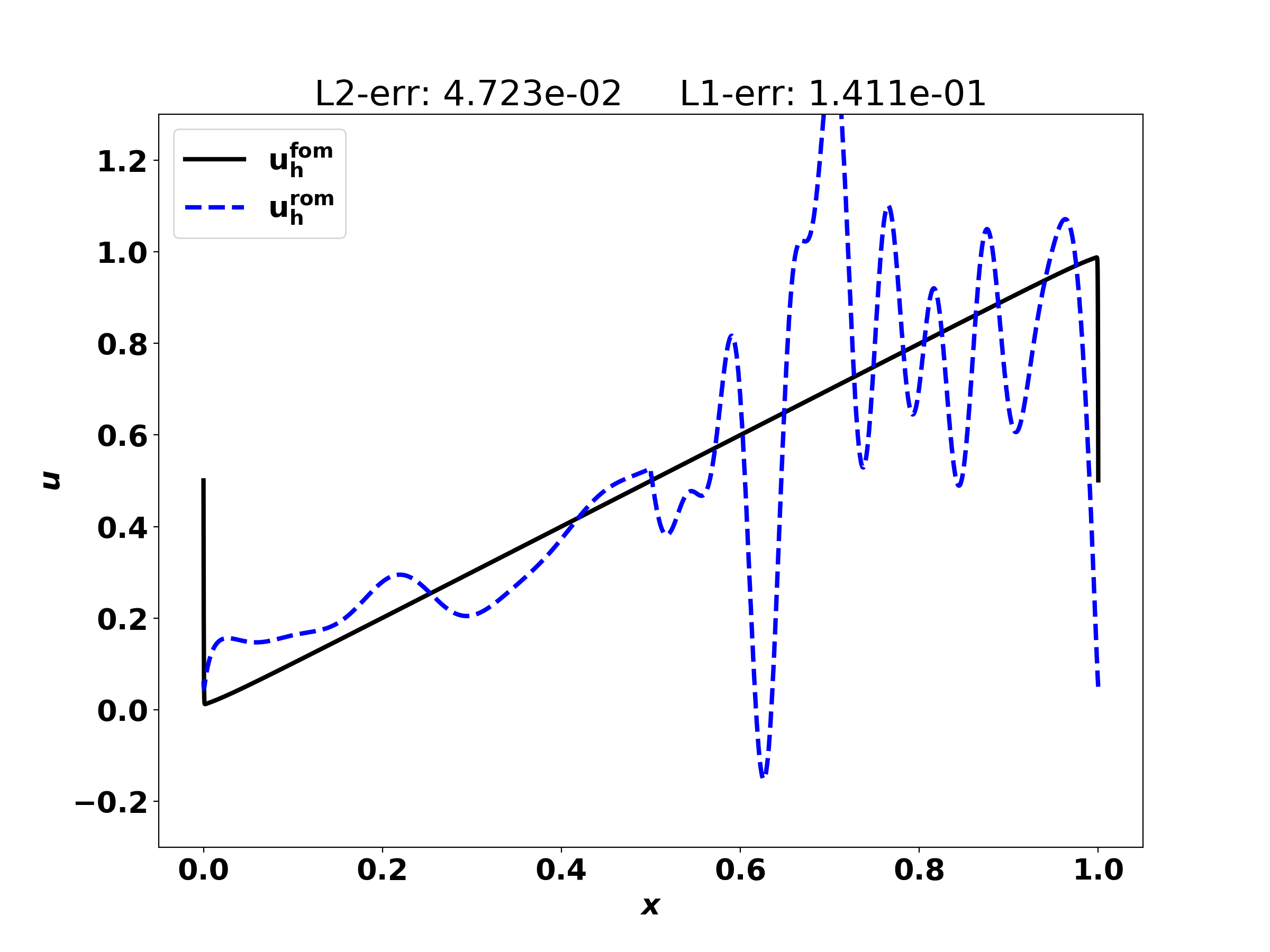}
  \caption{
Example 1.
    Numerical solution $u_h^{rom}$ for the
 POD-DG  model along with 
full order model solution $u_h^{fom}$. 
 Top: $t=0$. Middle: $t=0.5$. Bottom: $t=1$.
 Left: slightly resolved case $k=2$. 
 Right: fully resolved case $k=6$.
  }
\label{fig:ex1-soln}
\end{figure}

To build the POD model, we use $r=20$ basis functions, which capture
about $97.87\%$   of the total energy for both cases.
Numerical results for 
the plain POD-DG scheme \eqref{full-pod-online}, along with the computed
$L^2$- and $L^1$-errors for $u_h^{rom}-u_h^{fom}$ at time $t=0,0.5,1$, where 
$u_h^{fom}$ is the solution to the full order model \eqref{full-b}, and
$u_h^{rom}$ is the solution to the POD-DG model \eqref{full-pod-online},
are shown in Figure~\ref{fig:ex1-soln}.
It is clearly seen that the plain POD-DG model produces very oscillatory
results, with the associated error for $t=0.5$ and $t=1$ being an order of
magnitude larger than the initial projection error at $t=0$.

\begin{figure}[!ht]
\centering
\includegraphics[width=.45\textwidth]{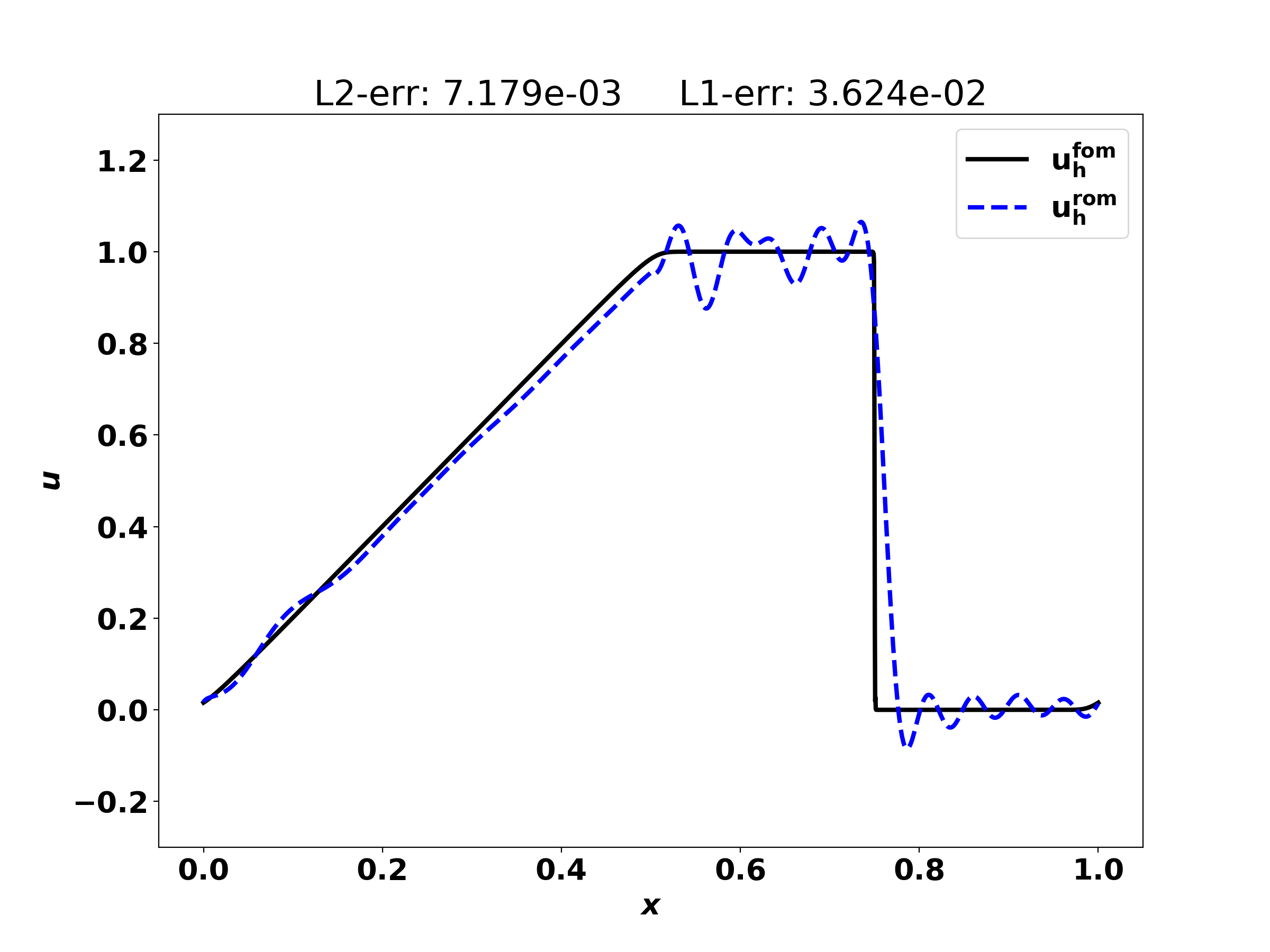}
\includegraphics[width=.45\textwidth]{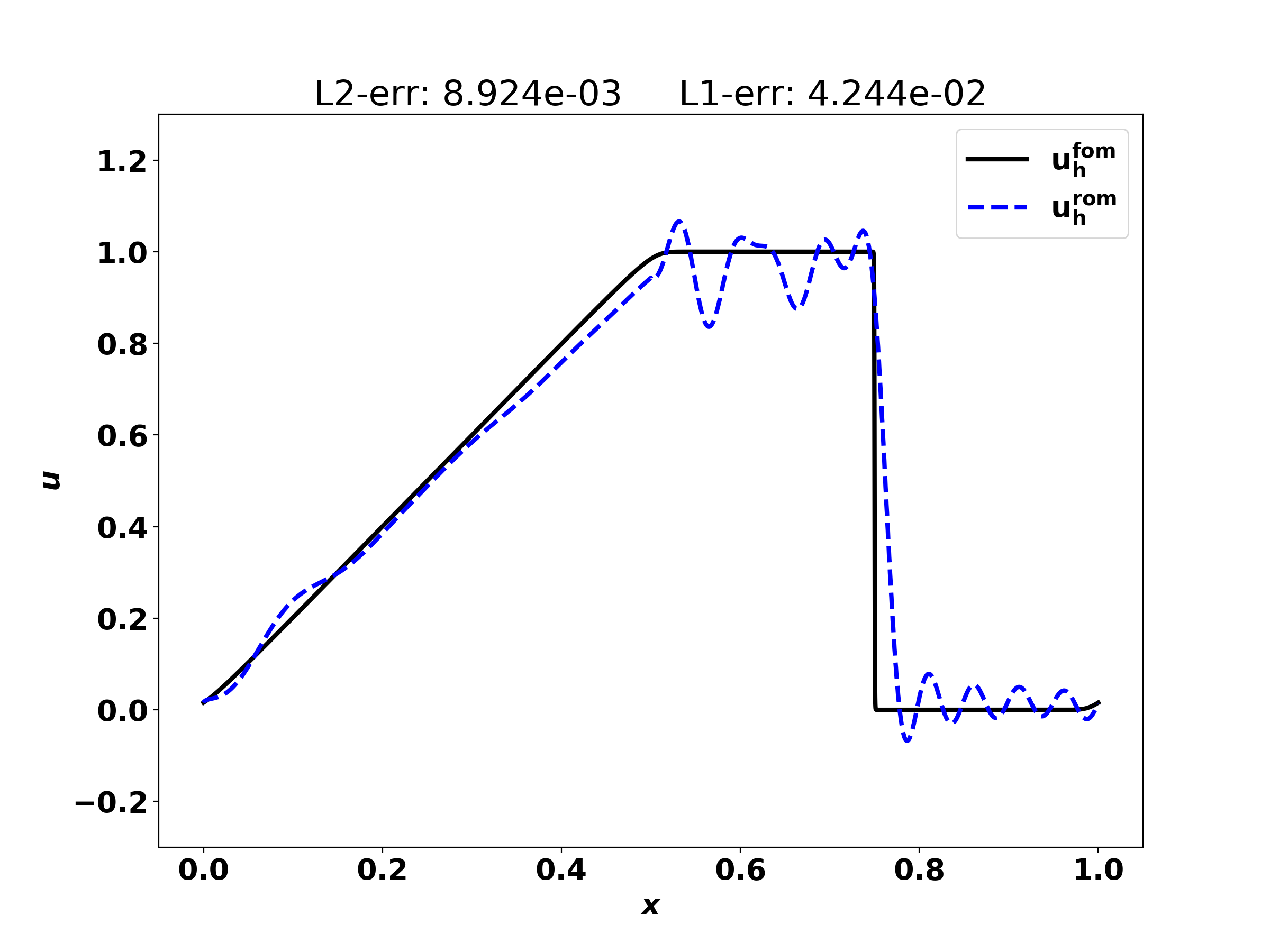}
\includegraphics[width=.45\textwidth]{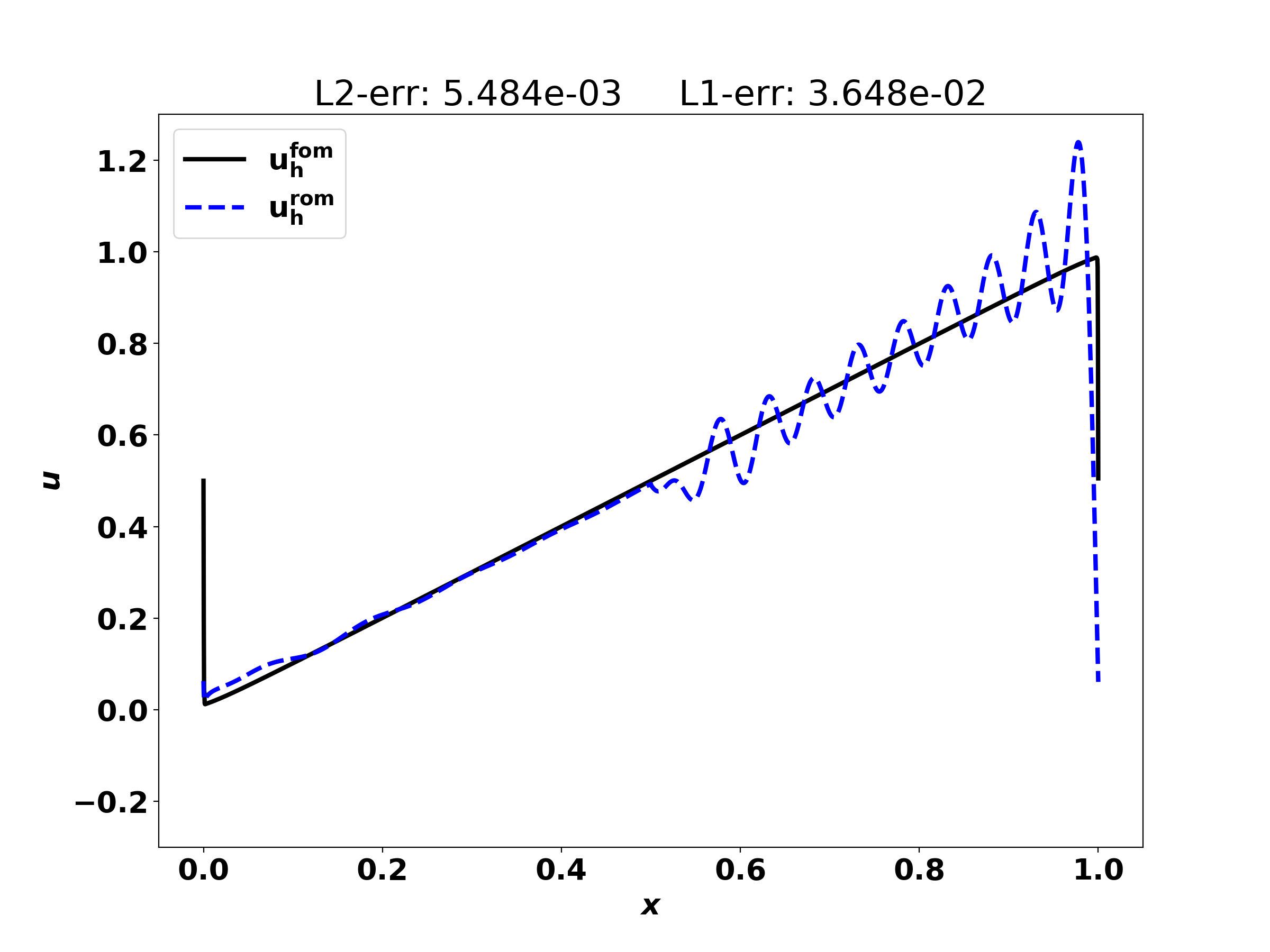}
\includegraphics[width=.45\textwidth]{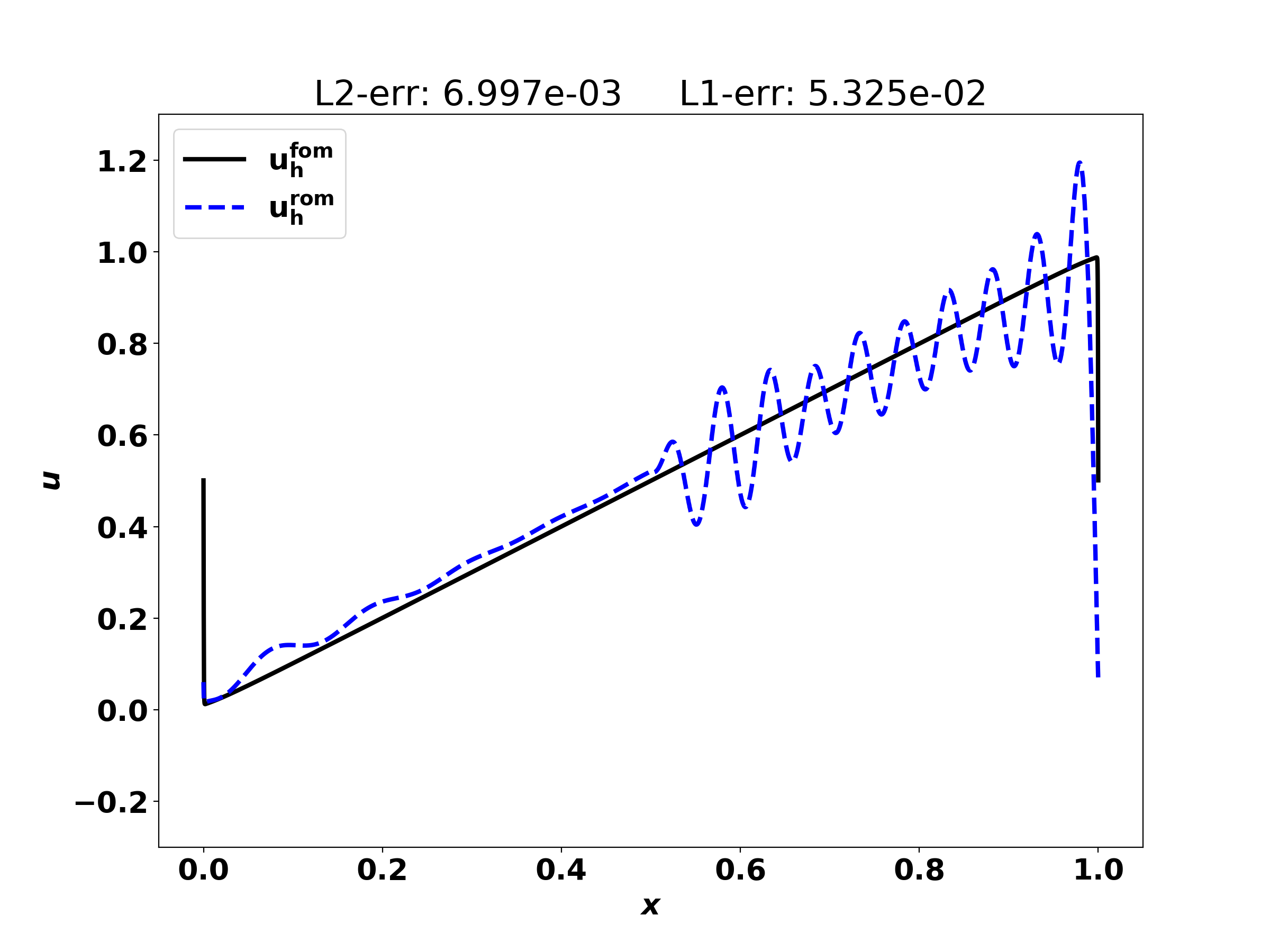}
  \caption{
Example 1.
    Numerical solution $u_h^{rom}$ for the
 POD-DG-C  model along with 
    full order model solution $u_h^{fom}$. 
 Top: $t=0.5$. Bottom: $t=1$.
 Left: slightly resolved case $k=2$, $c_1=10^4$. 
 Right: fully resolved case $k=6$, $c_1=2\times 10^8$.
  }
\label{fig:ex1-soln1}
\end{figure}

Next, we consider the POD-DG closure model \eqref{closure} with only convective
stabilization ($c_1>0, c_2=0$). We refer to the
resulting model as the POD-DG-C model.
We tune the parameter $c_1=10^4$ for $k=2$, and $c_1=2\times
10^8$ for $k=6$ to produce satisfactory results. Note that for the
classical upwinding DG
scheme \eqref{semi-b}, the parameter $c_1$ corresponds to the magnitude of
the solution
$|\avg{u_h}|\approx 1$, which is too small for the POD-DG model to suppress
numerical 
oscillation. We don't have a physical interpretation for the parameter
$c_1$, but argue that our global POD DG basis functions are very smooth across
element boundaries (which is especially true for the fully resolved case
$k=6$), and one needs to have a large weighting coefficient $c_1$ to 
make the convective stabilization term effective.
The associated numerical
results at $t=0.5$ and $t=1$ are shown in Figure~\ref{fig:ex1-soln1}.
Significant improvement over the results of the plain POD-DG model can be
clearly observed. 
We also found that the errors at $t=0.5$ and $t=1$ for the POD-DG-C
closure model are of similar magnitude to the POD projection
error at $t=0$ in Figure~\ref{fig:ex1-soln}.
However, the POD solution is still oscillatory behind the shock.

\begin{figure}[!ht]
\centering
\includegraphics[width=.45\textwidth]{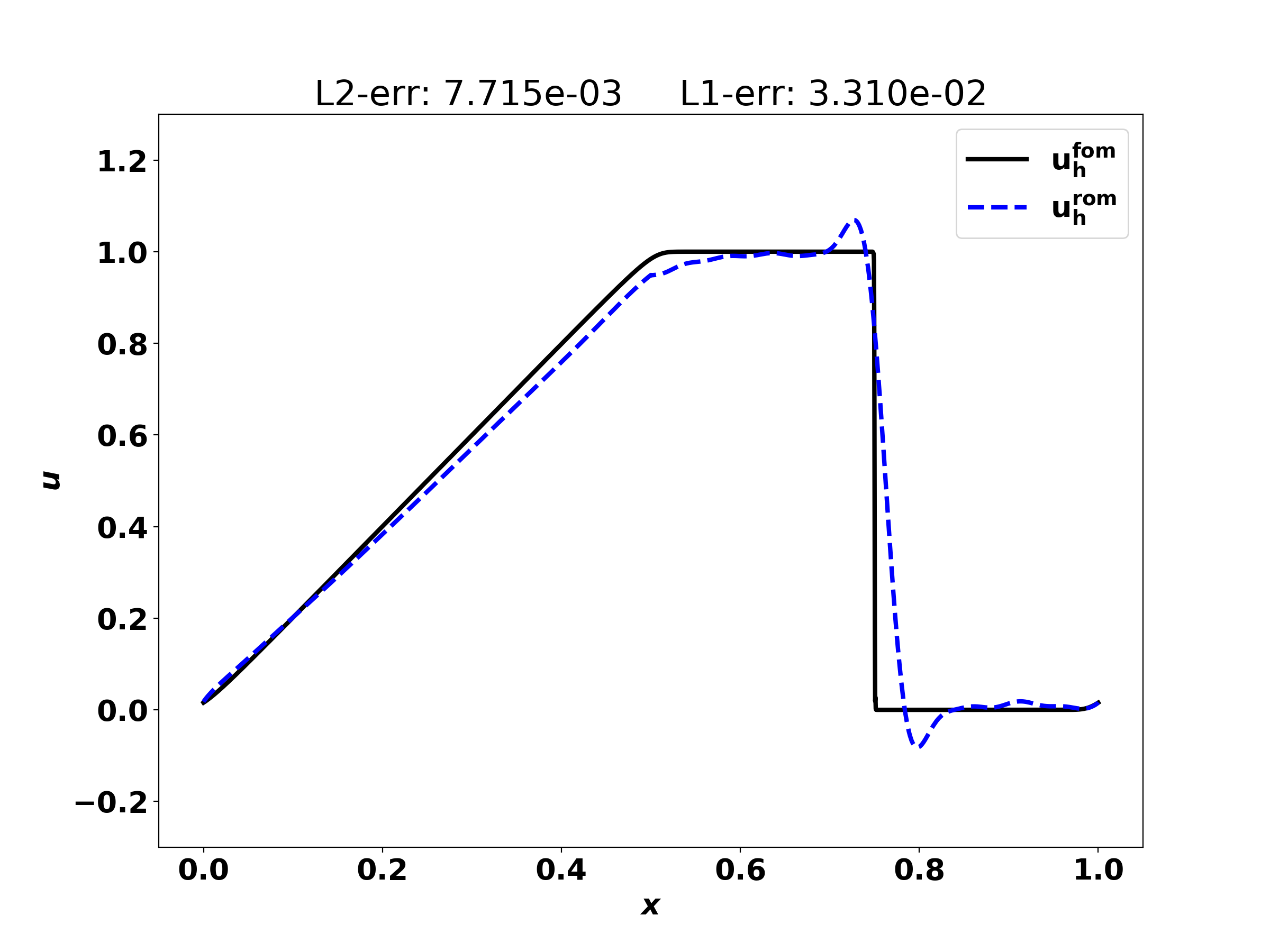}
\includegraphics[width=.45\textwidth]{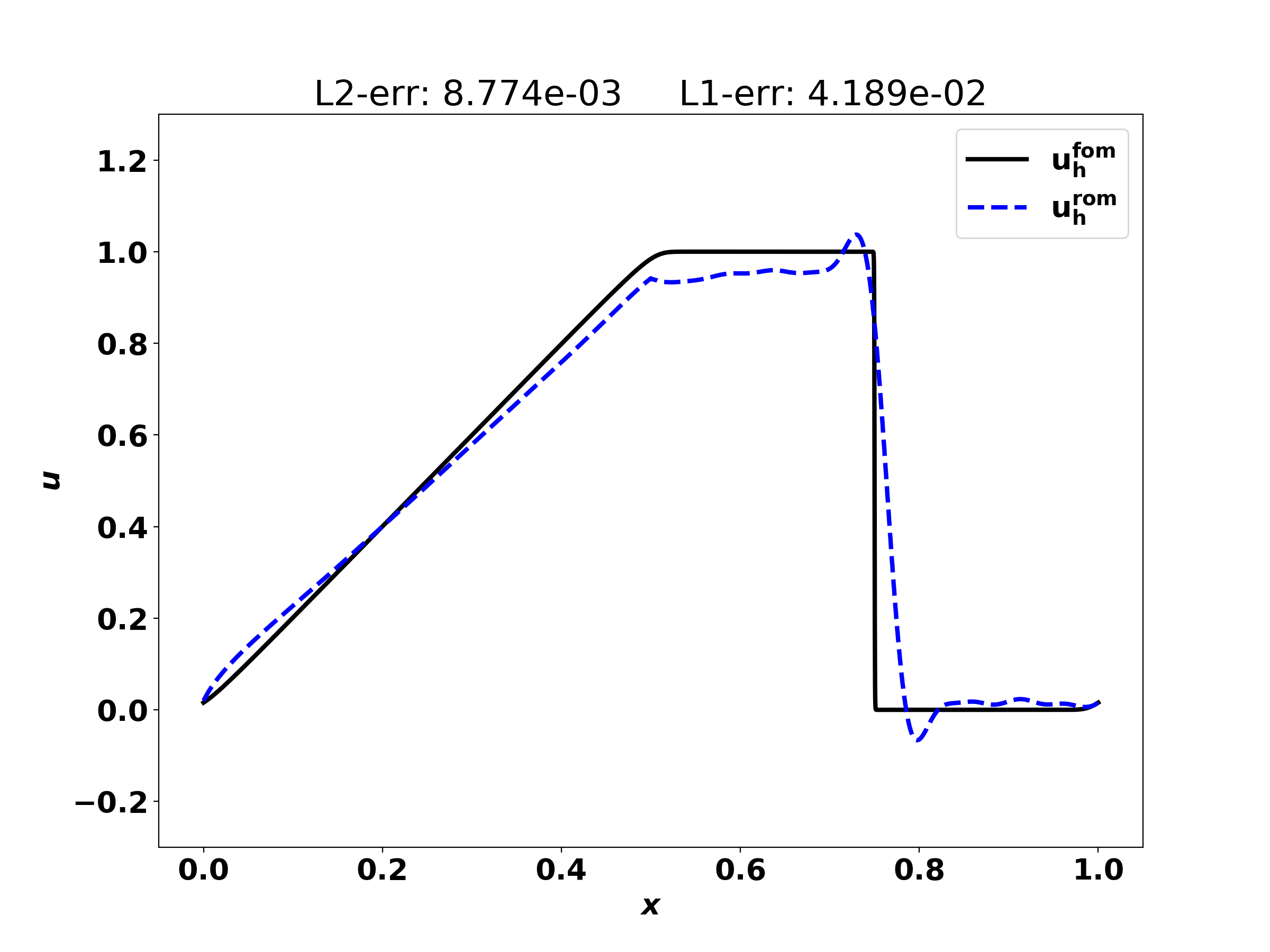}
\includegraphics[width=.45\textwidth]{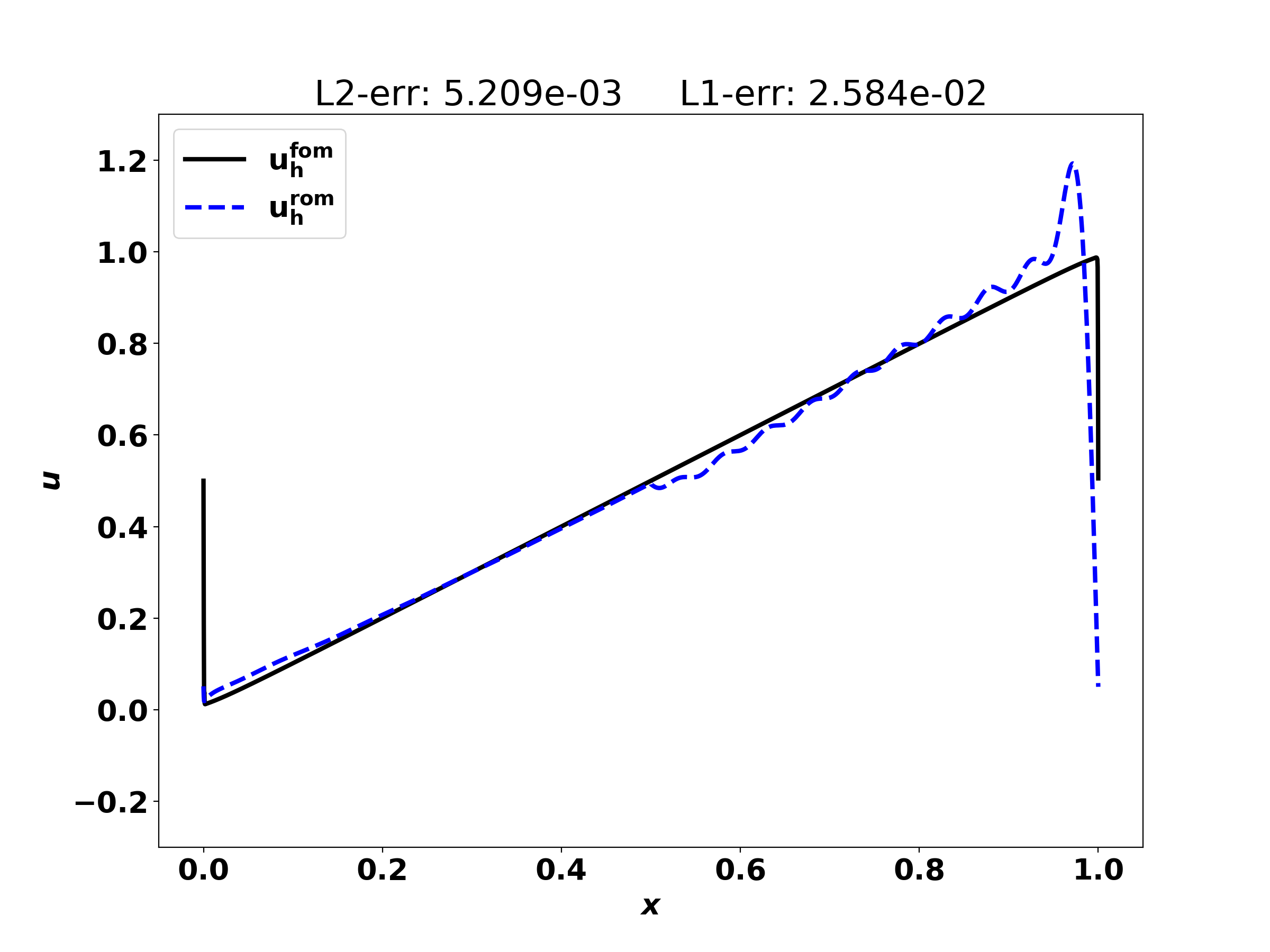}
\includegraphics[width=.45\textwidth]{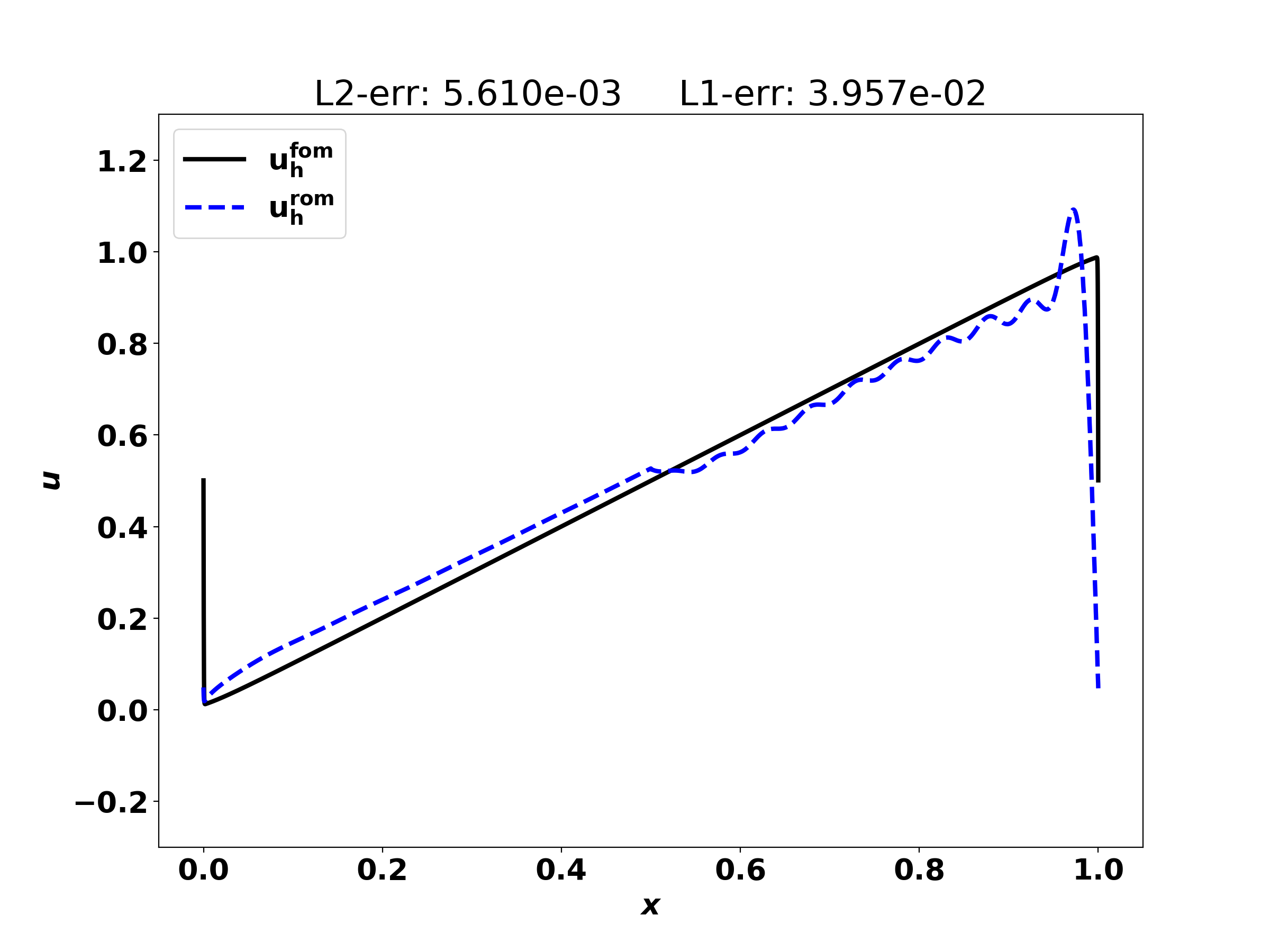}
  \caption{
Example 1.  Numerical solution $u_h^{rom}$ for the
 POD-DG-CD model along with 
    full order model solution $u_h^{fom}$. 
 Top: $t=0.5$. Bottom: $t=1$.
 Left: slightly resolved case $k=2$, $c_1=10^4, c_2=0.01$. 
 Right: fully resolved case $k=6$, $c_1=2\times 10^8, c_2=0.01$.
  }
\label{fig:ex1-soln2}
\end{figure}

Furthermore, we consider the POD-DG closure model \eqref{closure} with both convective and
diffusive stabilizations, which is referred to the POD-DG-CD model. 
We use the same parameter $c_1$ as the POD-DG-C model, i.e. $c_1=10^4$ for $k=2$, and $c_1=2\times
10^8$ for $k=6$; and set $c_2=0.01$.
The associated numerical
results at $t=0.5$ and $t=1$ are shown in Figure~\ref{fig:ex1-soln2}.
We observe that the errors at $t=0.5$ and $t=1$ for the POD-DG-CD model is
similar and
slightly smaller than those for the POD-DG-C model, and the post-shock
oscillations are also diminished.

\begin{figure}[!ht]
\centering
\includegraphics[width=.45\textwidth]{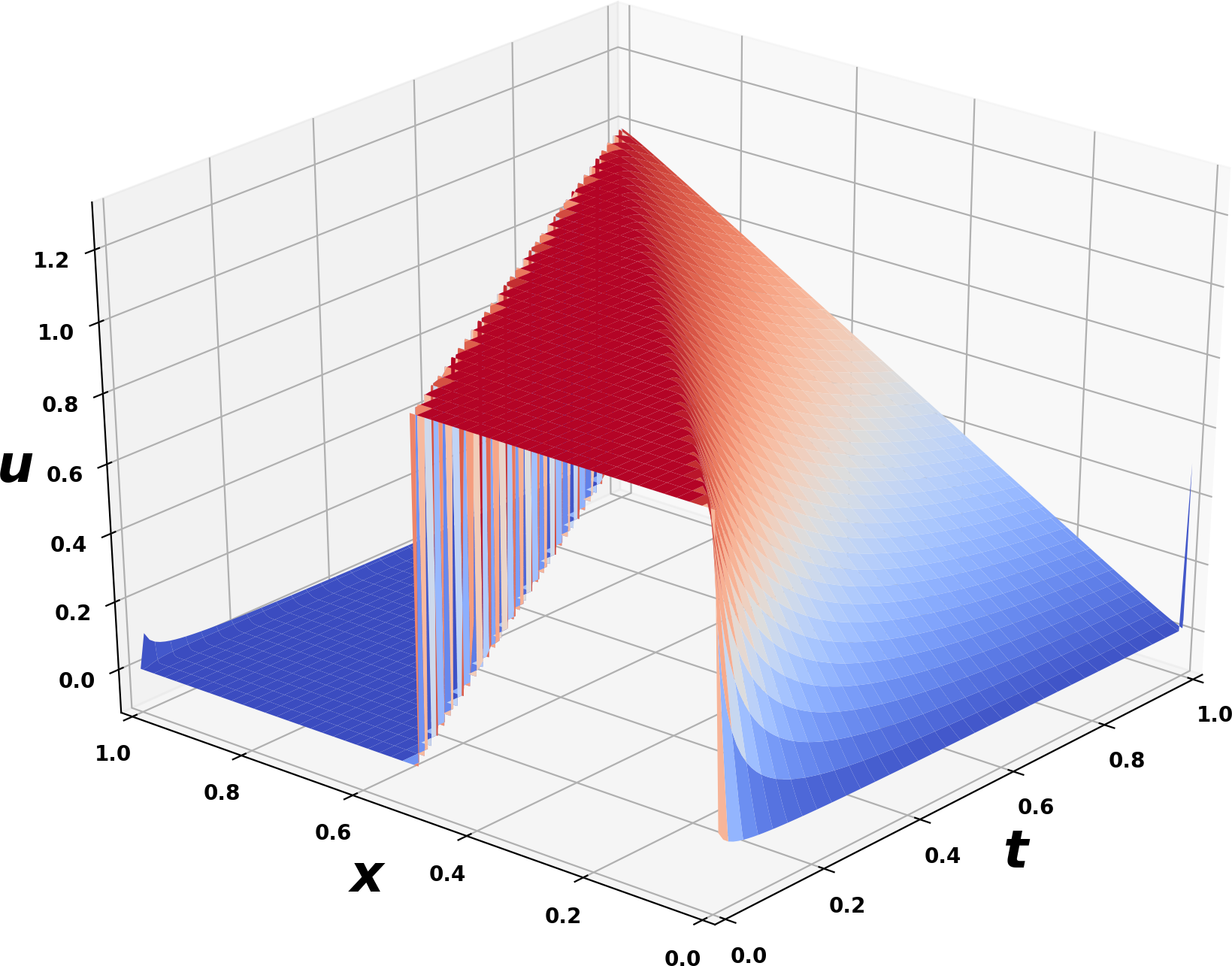}
\includegraphics[width=.45\textwidth]{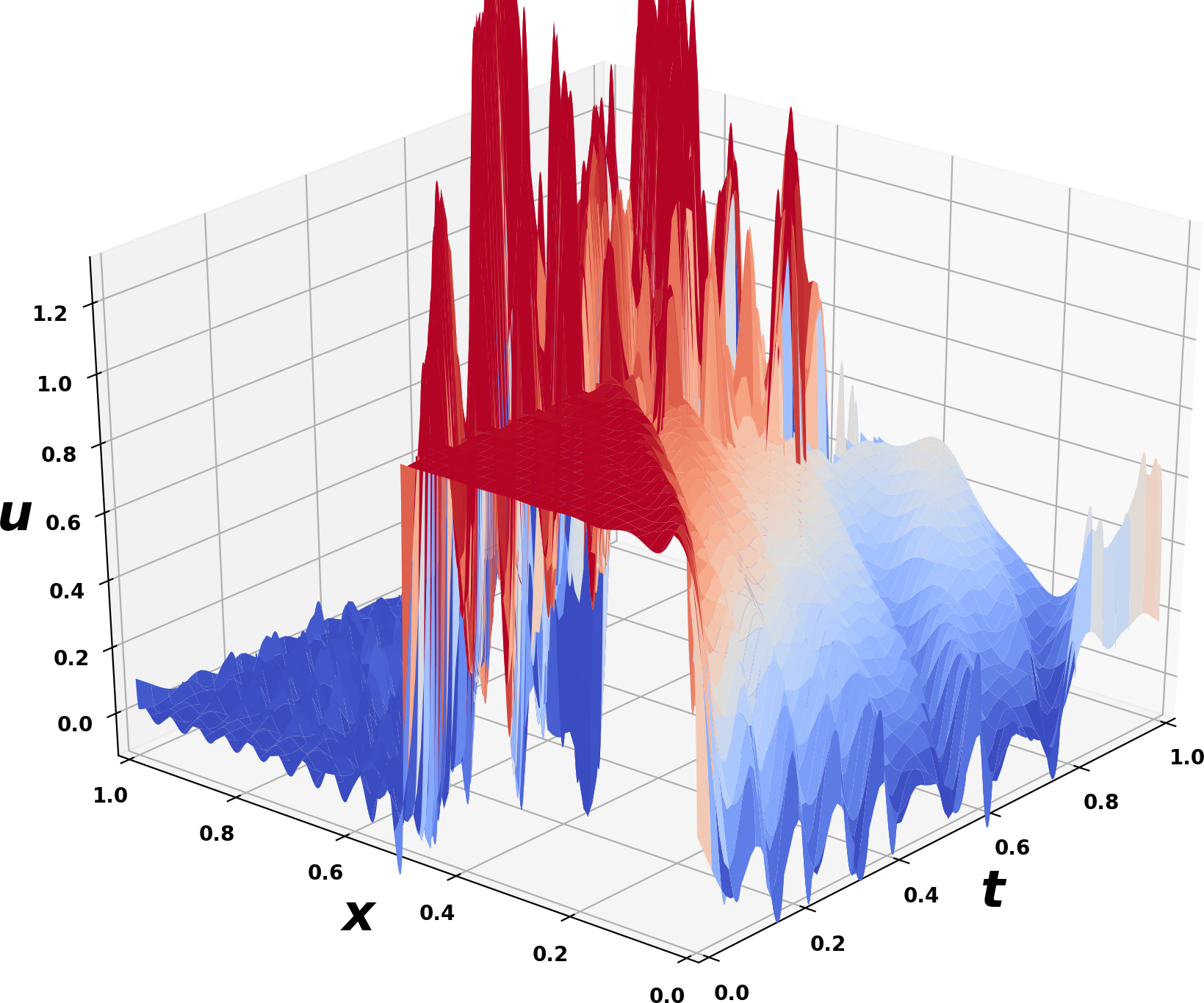}
\includegraphics[width=.45\textwidth]{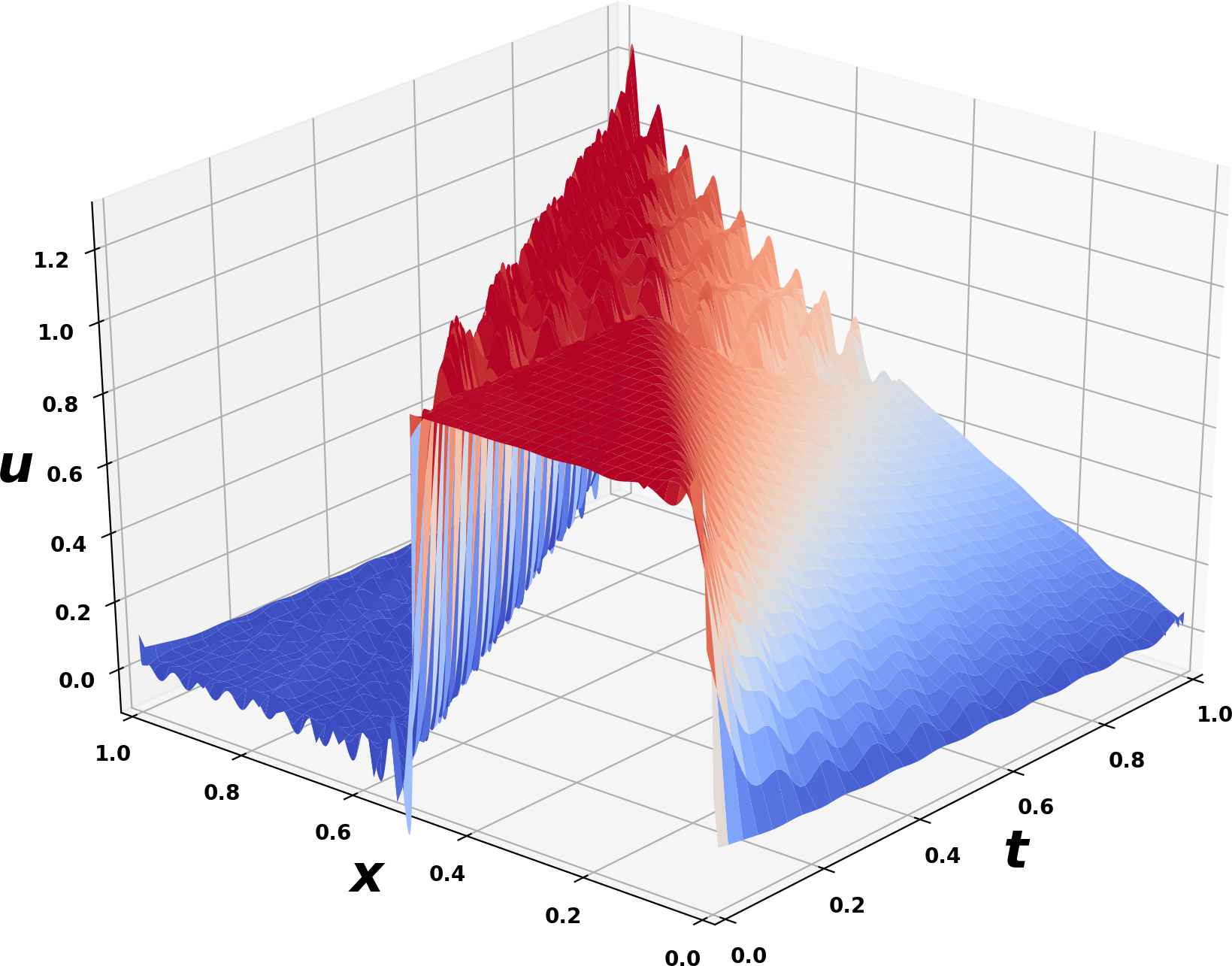}
\includegraphics[width=.45\textwidth]{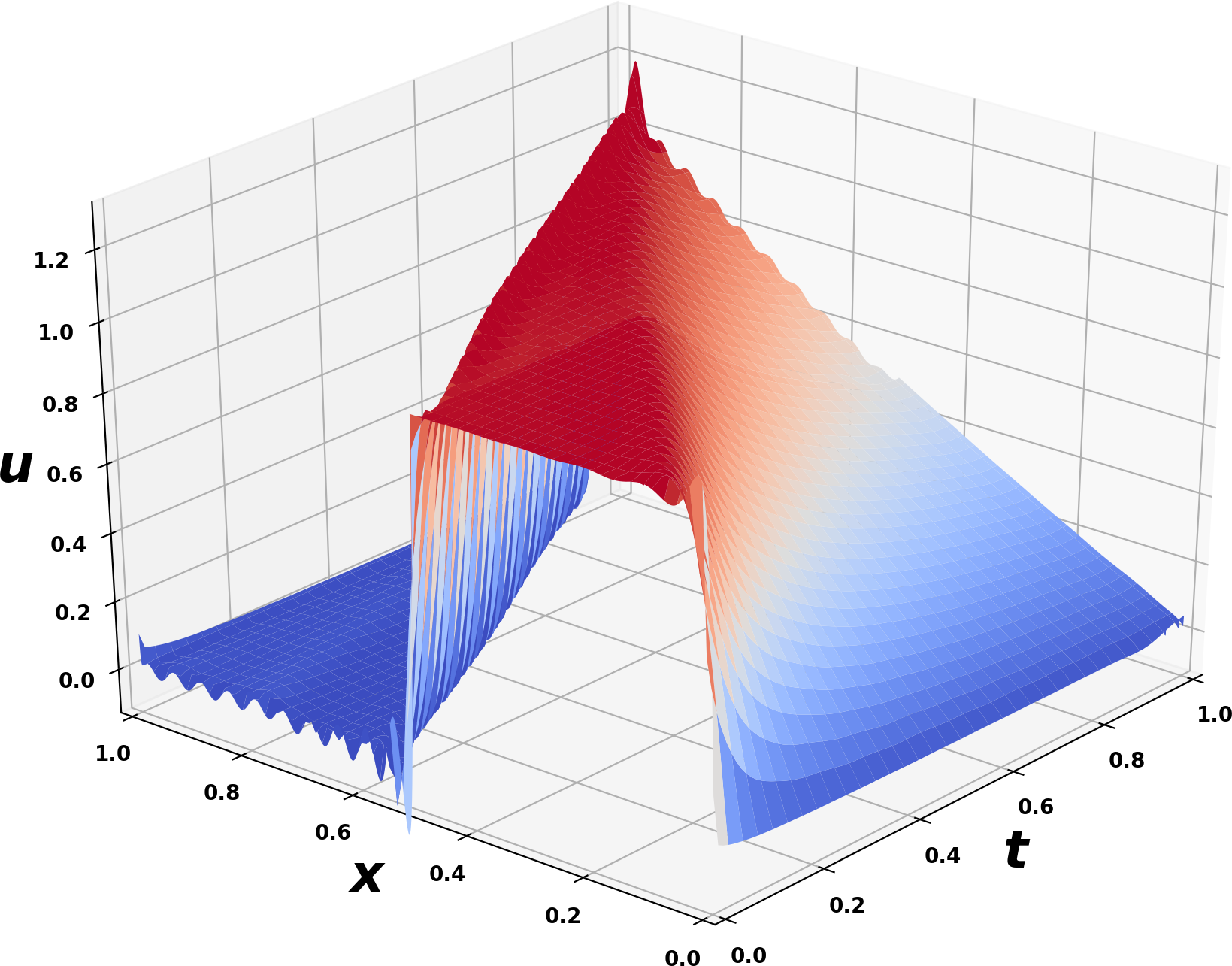}
  \caption{
Example 1.  Time evolution of numerical solutions.
  Top left: full order model. Top right: POD-DG. Bottom left: POD-DG-C.
  Bottom right: POD-DG-CD. 
  Slightly resolved case $k=2$.
  }
\label{fig:ex1-time1}
\end{figure}

Finally, the time evolution of the three models along with the full order model are
presented in Figure~\ref{fig:ex1-time1} for $k=2$. 
The results for $k=6$ are similar and are omitted to save space.


\subsection{Example 2. Burger's equation: smooth initial condition}
\begin{figure}[!ht]
\centering
\includegraphics[width=.32\textwidth]{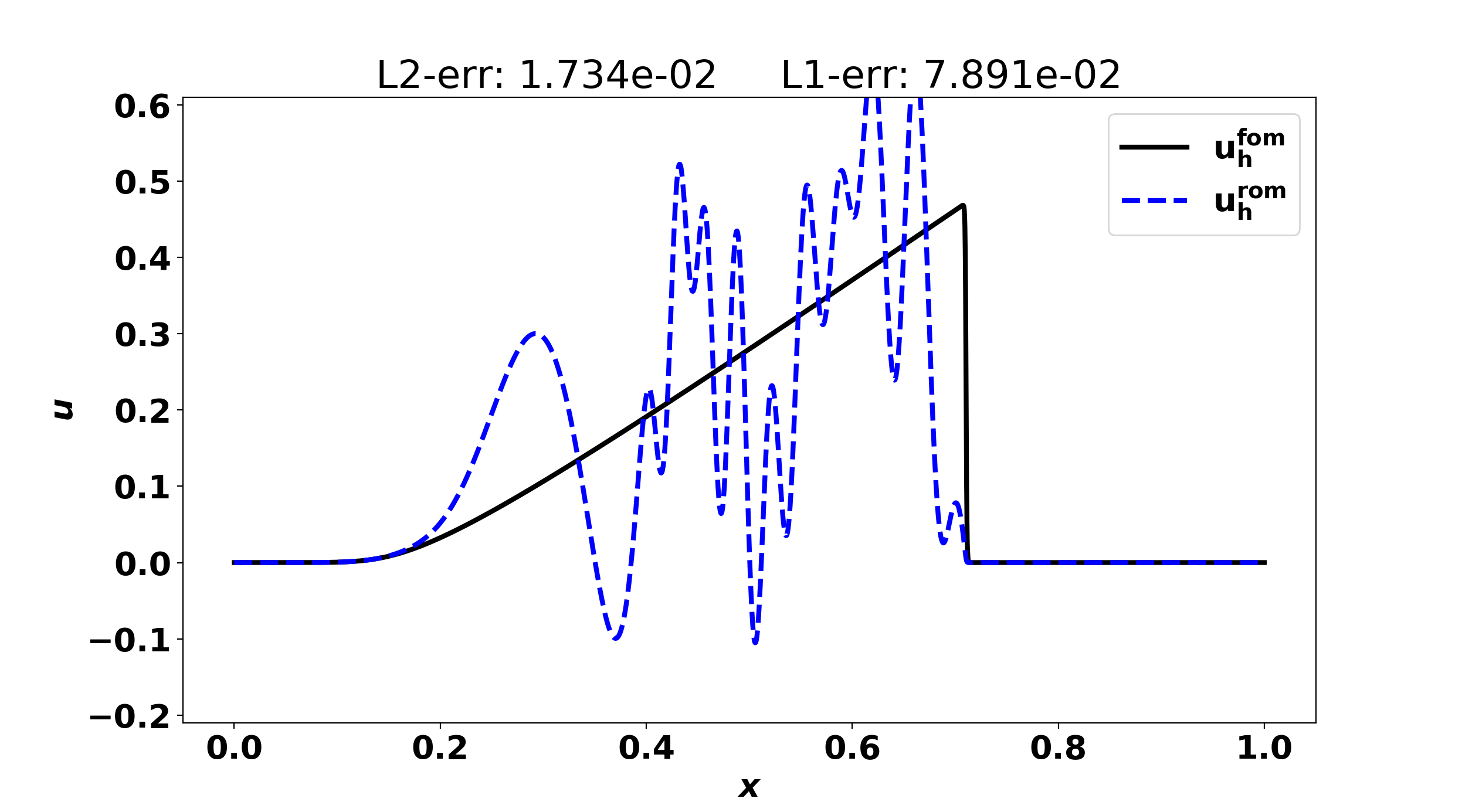}
\includegraphics[width=.32\textwidth]{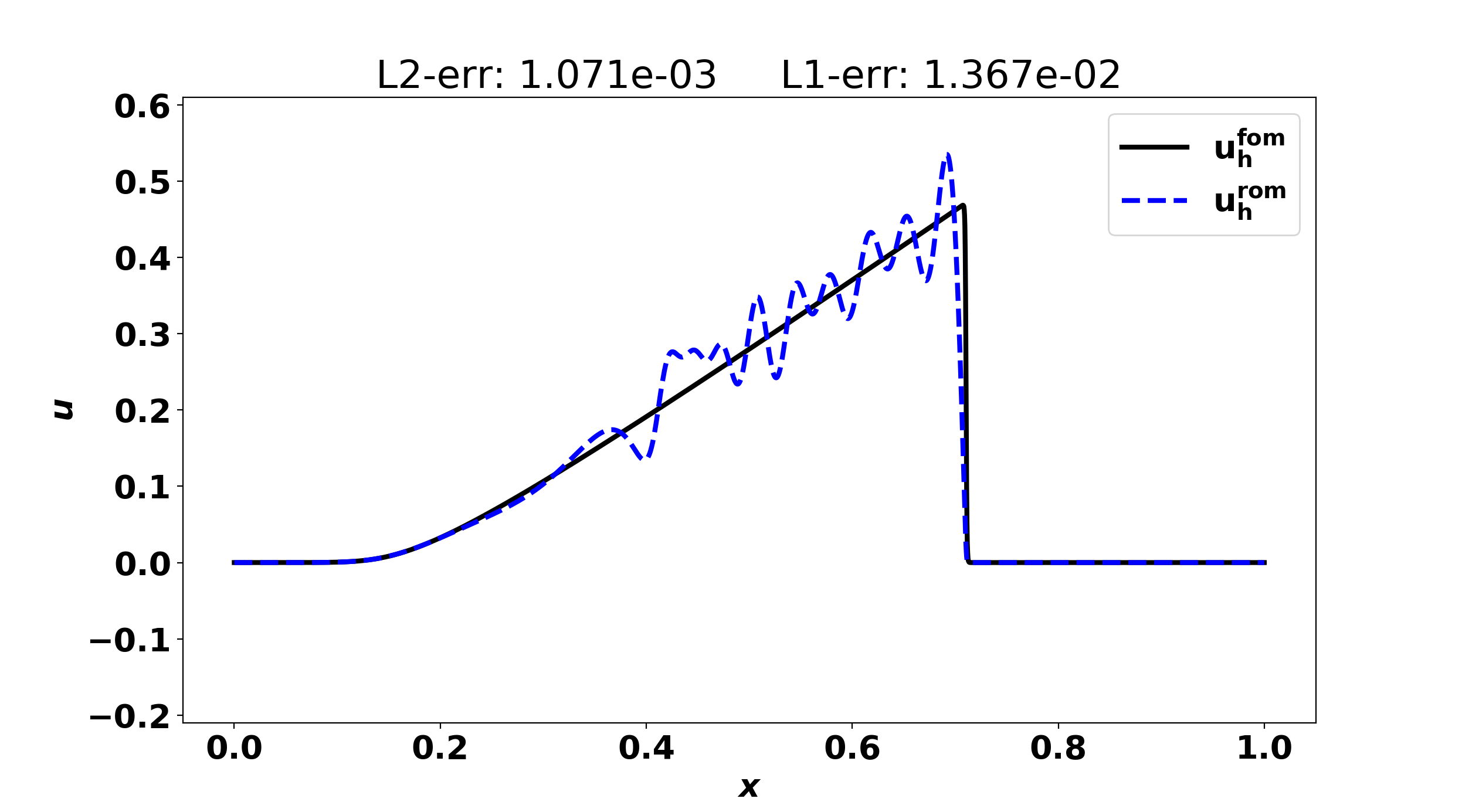}
\includegraphics[width=.32\textwidth]{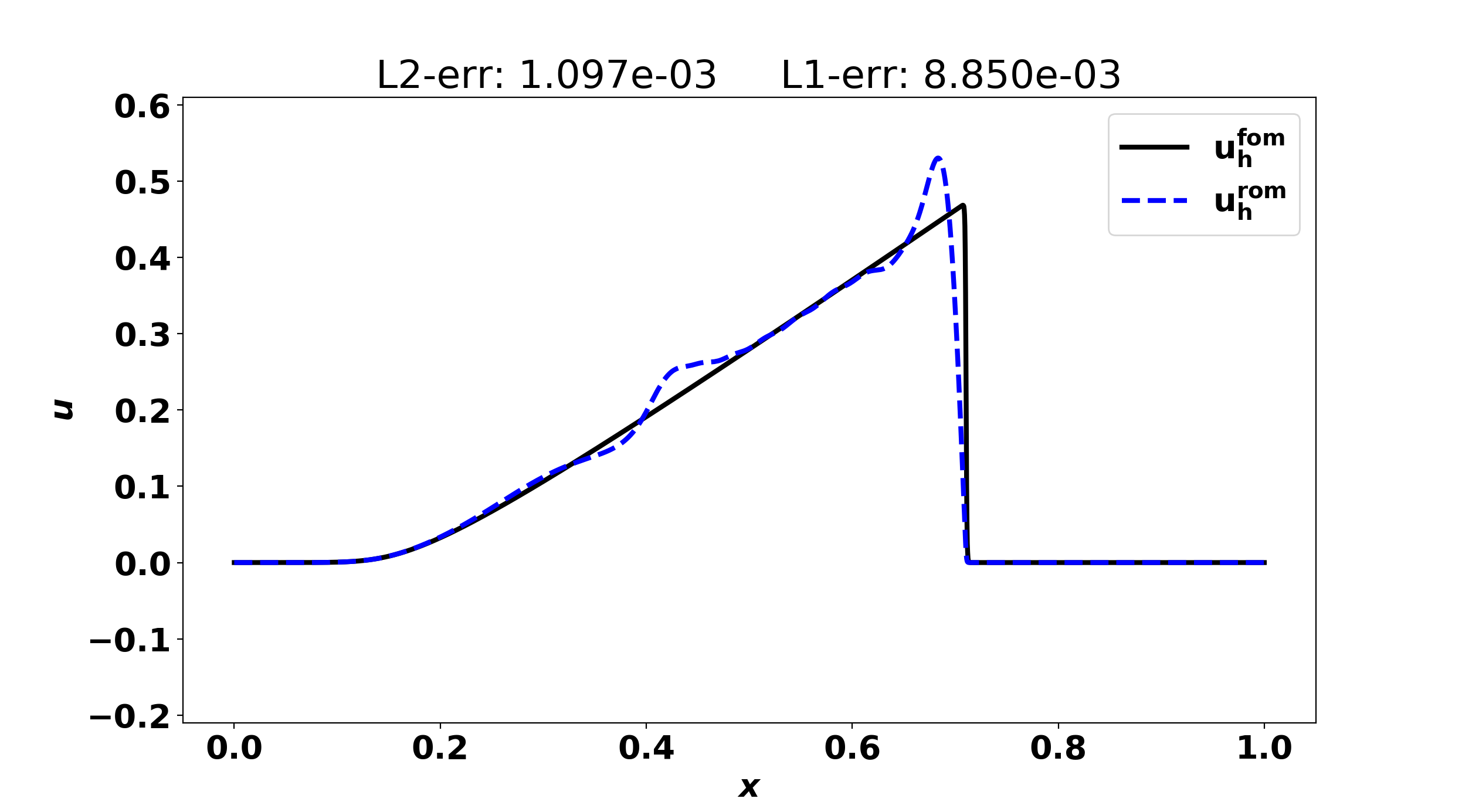}
  \caption{
Example 2. Numerical solution at final time $t=1$. 
 Left: POD-DG model. Middle: POD-DG-C model 
 with $c_1=10^4,c_2=0$. 
 Right: POD-DG-CD model
 with $c_1=10^4,c_2=0.01$.
 20 POD bases are used.
 Slightly resolved case $k=2$, $h=10^{-4}$.
  }
\label{fig:ex2-soln}
\end{figure}

We consider the same problem as Example 1, but with the following smooth initial
condition:
  \[
    u(0) = \exp\left({-200(x-0.3)^2}
    \right).
  \] 
 Very similar results as those for Example 1 are observed. 
 In particular, we need to take  $c_1=10^4$ for the case $k=2$,
 and $c_2=10^8$ for the case $k=6$ to make the POD-DG-C model 
 produce satisfactory results, and use the POD-DG-CD model with $c_2=0.01$ 
 to further improve the results. 
 We present in Figure~\ref{fig:ex2-soln} the numerical solution of different
 POD-DG models at final time $t=1$ for the slightly resolved case $k=2$.
 It is again clear that the POD-DG-C model produce better results than the
 plain POD-DG model, and the 
 POD-DG-CD model further improves the results of POD-DG-C model by suppressing 
 post-shock oscillations.

\subsection{Example 3. Navier-Stokes: 
2D flow past a cylinder, $Re=100$}
\begin{figure}[ht!]
 \includegraphics[width=.75\textwidth]{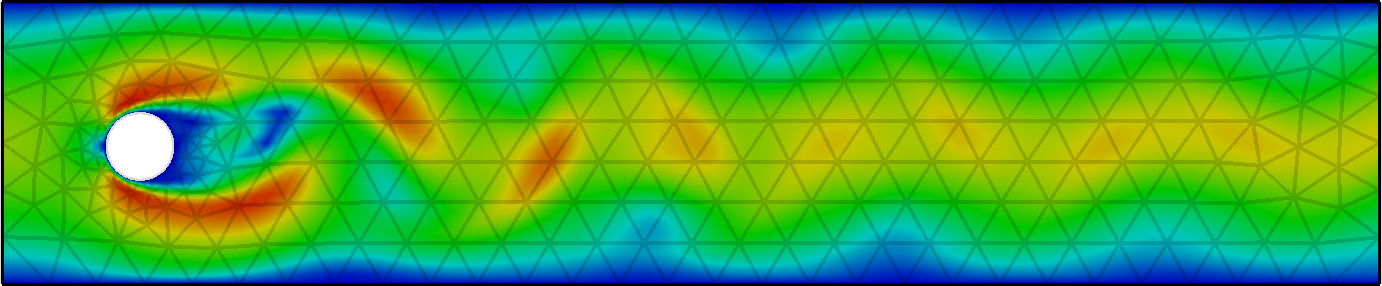}
 \vspace{1ex}
  \caption{Example 3: the initial velocity field
 (color corresponding to velocity magnitude
 $\|\bld u_h\|_2$ from $0$ to $2.17$).}
 \label{fig:obs}
\end{figure}

We consider the classical flow past a cylinder benchmark problem 
\cite{LD96}. 
The domain is a rectangular channel with an almost vertically centered circular obstacle, 
c.f. Fig.~\ref{fig:obs},
\[
 \Omega:=[0,2.2]\times [0,0.41]\backslash \{\|(x,y)-(0.2,0.2)\|_2\le 0.05\}.
\]
The boundary is decomposed into $\Gamma_{{in}}:=\{x=0\}$, the inflow boundary, 
$\Gamma_{{out}}:=\{x=2.2\}$, the outflow boundary, and 
$\Gamma_{{wall}}:=\partial\Omega\backslash(\Gamma_{{in}}\cup \Gamma_{{out}})$,
the wall boundary. On $\Gamma_{{out}}$ we prescribe 
natural boundary conditions
$(-\nu\nabla\bld u + p I) \bld n = 0$, on 
$\Gamma_{{wall}}$ homogeneous Dirichlet boundary conditions for the velocity (no-slip) 
and on $\Gamma_{{in}}$  the
inflow Dirichlet boundary conditions
\[\bld u(0, y, t) = 6\bar u\, y(0.41 - y)/0.41^2 \cdot (1, 0),
\]
with $\bar u=1$ the average inflow velocity.
The viscosity is taken to be $\nu = 10^{-3}$, hence Reynolds number $\mathrm{Re}=\bar u D/\nu = 100$,
where $D=0.1$ is the disc diameter.

For this Reynolds number, the flow turns into a time-periodic behaviour
with a vortex shedding behind the cylinder. 
For the FOM, we consider the scheme \eqref{full-n} 
with polynomial degree $k=3$ on a (curved) unstructured triangular mesh 
with
292 triangular elements, and take time step size $\Delta t = 0.001$.
A precomputed fully developed velocity profile is
used for the initial condition; see Fig. \ref{fig:obs} for the geometry, 
the mesh and the initial velocity field.

\begin{figure}[ht!]
 \includegraphics[width=.75\textwidth]{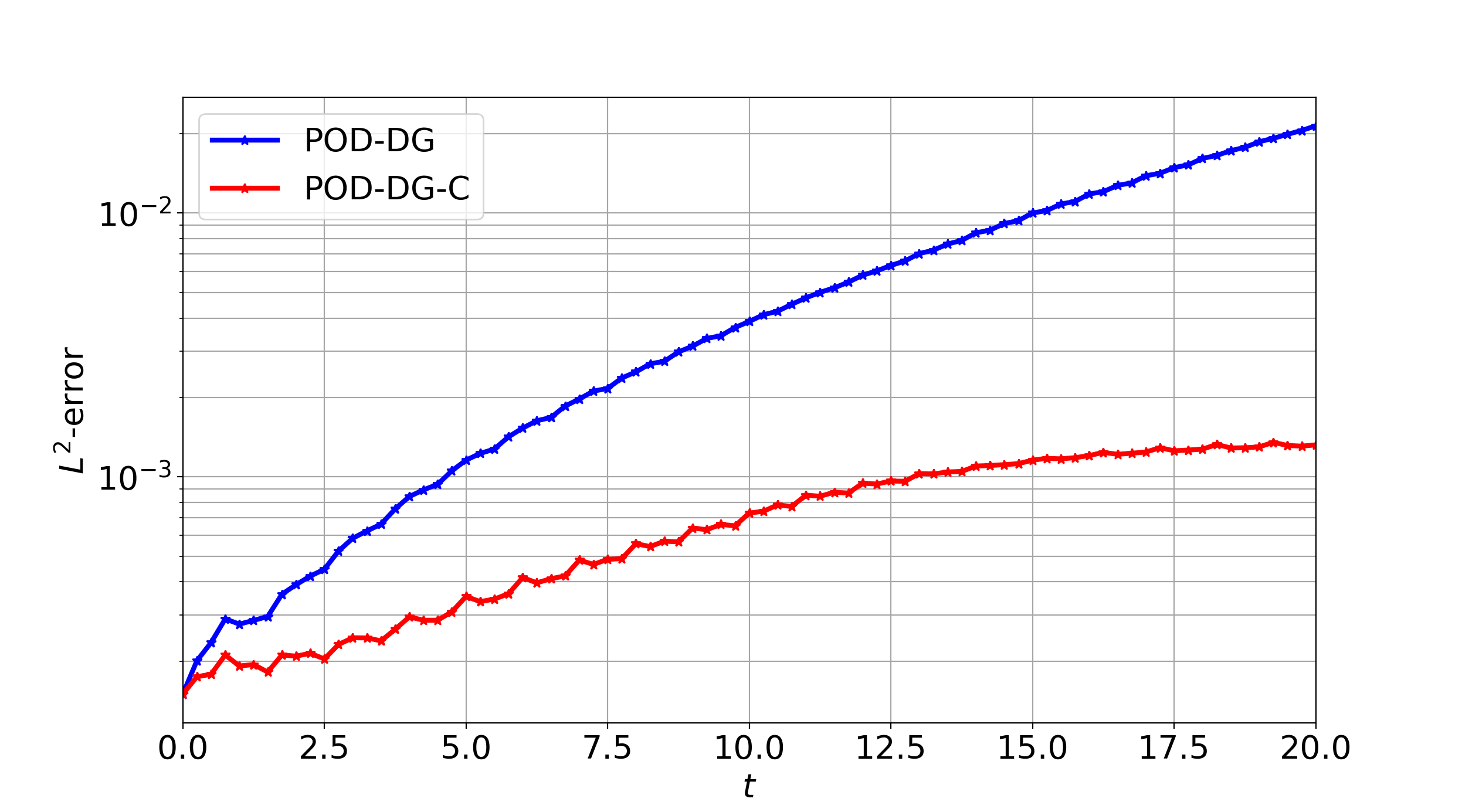}
  \caption{Example 3: time evolution of 
  the $L^2$-velocity error 
$\|\bld
u_h^{fom}-\bld u_h^{rom}\|$.} 
 \label{fig:ns-err}
\end{figure}

To build the POD bases, we 
collect $401$ snapshots in the time interval
$[0,2]$ taken at equidistant time instance. 
To build the POD model, we use 6 POD bases which capture about 
99.81\% of the total energy and run the simulation up to time 
$T=20$. We consider the plain POD-DG model and the POD-DG-C model with
$c_1=5$. 
The constant $c_1$ is tuned to yield relatively the
smallest $L^2$-error between FOM and ROM solutions at final time for
a range of choices. Taking $c_1$ too big or too small leads to less
accurate approximations. It is interesting to observe that this time 
$c_1$ is close to the maximum velocity magnitude $v_{\max}\approx 2.17$, 
which is very different to the scaling in the Burgers' equation cases in 
Examples 1-2.
Here, probably due to the relative small Reynolds number, we find
that adding extra diffusive stabilization in \eqref{closure} does not
improve the results. Hence, results for the POD-DG-CD model will not be
shown.
The time evolution of the $L^2$ velocity error 
$\|\bld
u_h^{fom}-\bld u_h^{rom}\|$ 
is plotted in Figure~\ref{fig:ns-err}.
We observe that the error for the POD-DG-C model is an order of
magnitude smaller than that for the plain POD-DG model at time $t=20$.

\begin{figure}[ht!]
 \includegraphics[width=.75\textwidth]{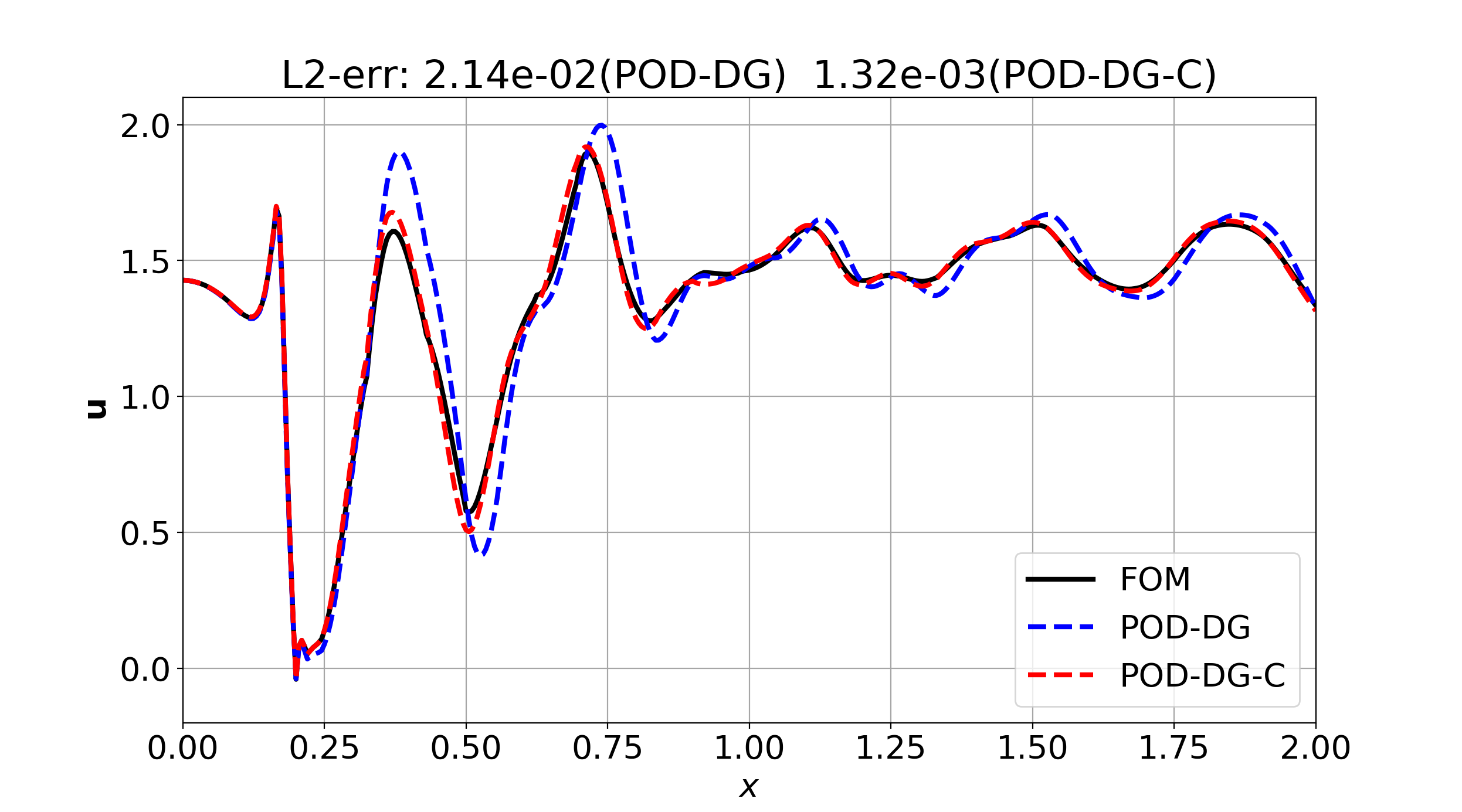}
  \caption{Example 3: x-component of velocity field along cut line 
  $y=0.25$ at time $t=20$.} 
 \label{fig:ns-u}
\end{figure}

\begin{figure}[ht!]
 \includegraphics[width=.45\textwidth]{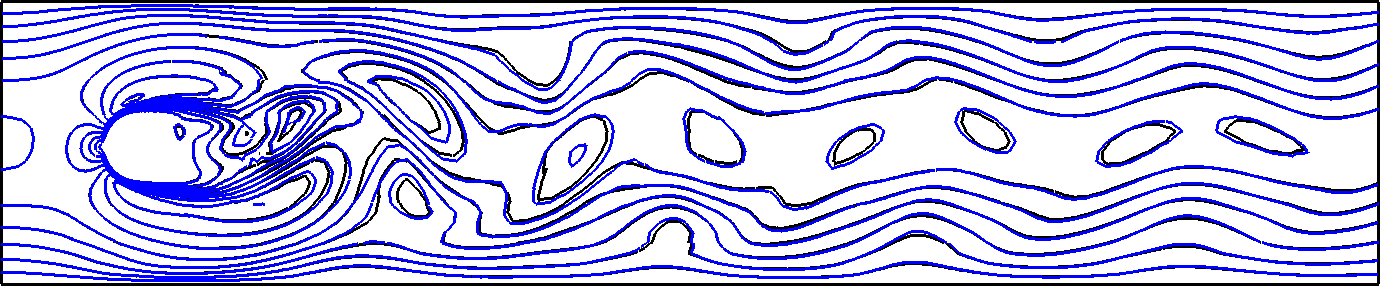}
 \includegraphics[width=.45\textwidth]{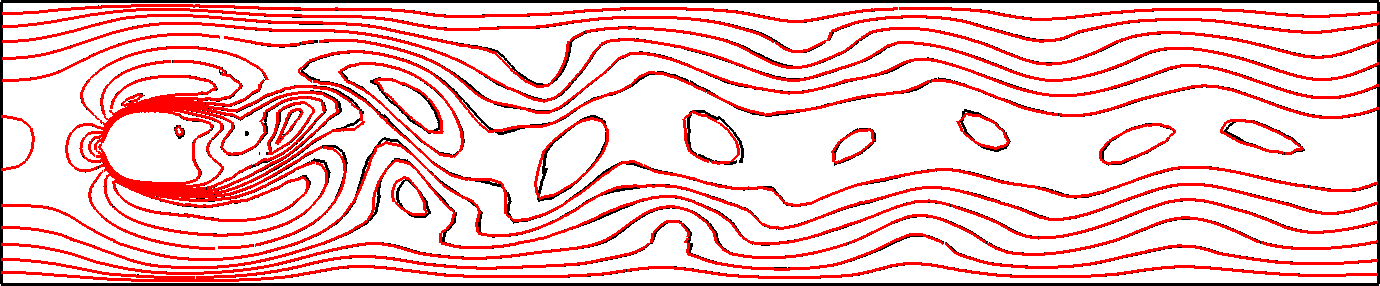}\\[.5ex]
 \includegraphics[width=.45\textwidth]{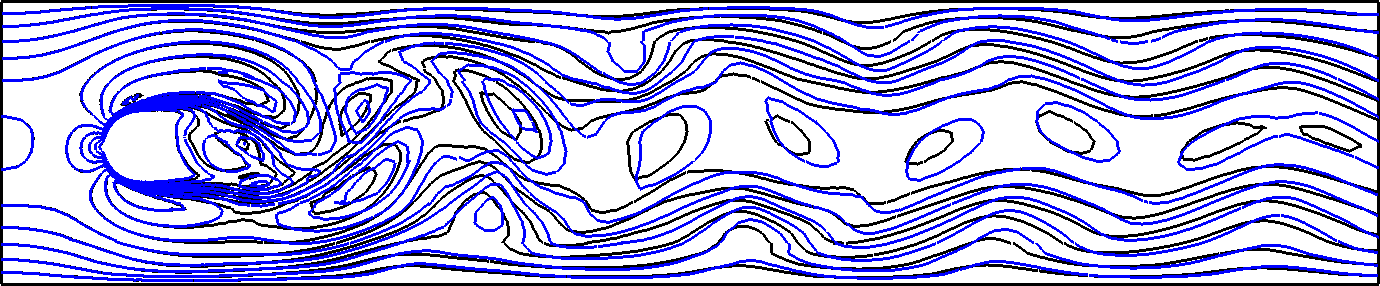}
 \includegraphics[width=.45\textwidth]{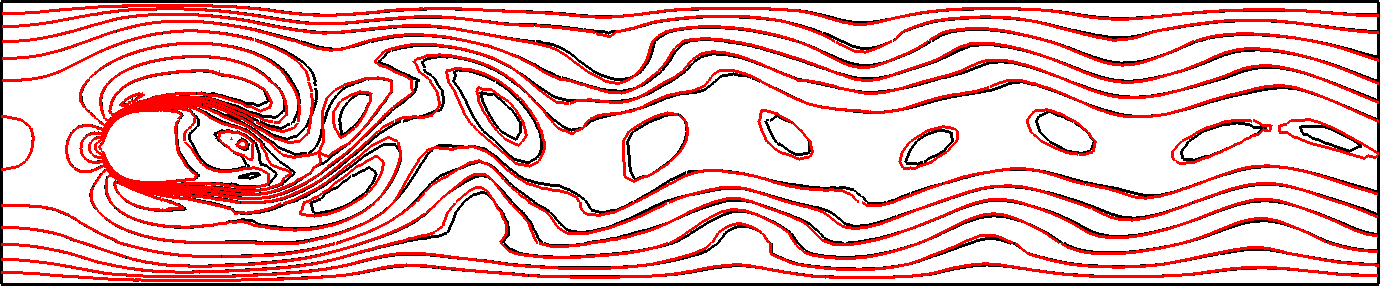}
 \caption{Example 3 (ROMs are 6-dimensional): velocity contour at time $t=5$ (top) and $t=20$
   (bottom). 10 equispaced contour lines from $0$ to $2.17$.
   Black: FOM. Bule: POD-DG. Red: POD-DG-C with $c_1=5$. 
} 
 \label{fig:ns-v}
\end{figure}

We plot 
the x-component of the velocity field along the 
cut line $y=0.25$ at time $t=20$ in Figure~\ref{fig:ns-u}. Clearly the result for the
POD-DG-C model is closer to FOM than that for the plain POD-DG model, which
produces a visible phase shift.
Finally, the velocity magnitude contour lines at time $t=5$ and $t=20$ for different models
are shown in Figure~\ref{fig:ns-v}. Here we observe that at time $t=5$,
both POD-DG and POD-DG-C models produce similar results as the FOM.
On the other hand, visible phase shift, especially behind the cylinder, is observed for the POD-DG model
(in blue) at time $20$, while the result for POD-DG-C (in red) is still in good
agreement with FOM.

\subsection{Example 4. Navier-Stokes: 
2D flow past a cylinder, $Re=500$}
\begin{figure}[ht!]
 \includegraphics[width=.75\textwidth]{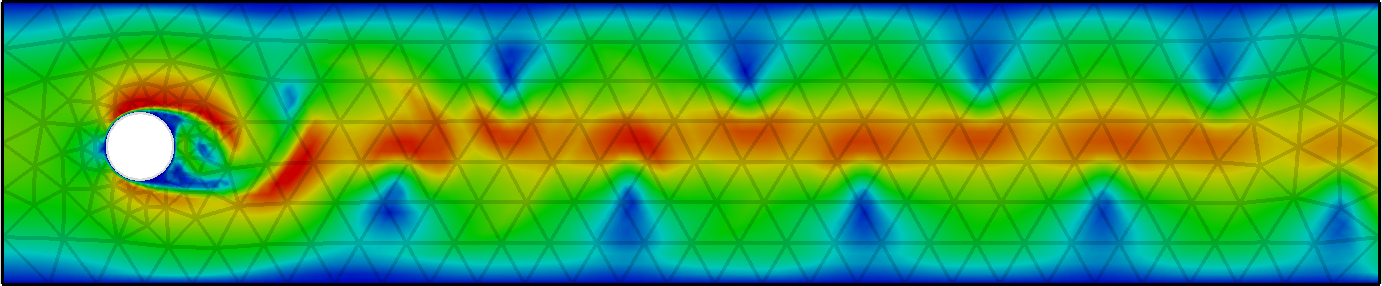}
 \vspace{1ex}
  \caption{Example 4: the initial velocity field 
 (color corresponding to velocity magnitude
 $\|\bld u_h\|_2$ from $0$ to $2.4$).}
 \label{fig:obs2}
\end{figure}

We consider the same problem as Example 3, but with a larger Reynolds number
$Re=500$. For the FOM, we consider the scheme \eqref{full-n} 
with polynomial degree $k=6$ on 
the mesh used in Example 3. The initial (fully developed)
velocity field is shown in  Figure~\ref{fig:obs2}.

\begin{figure}[ht!]
 \includegraphics[width=.75\textwidth]{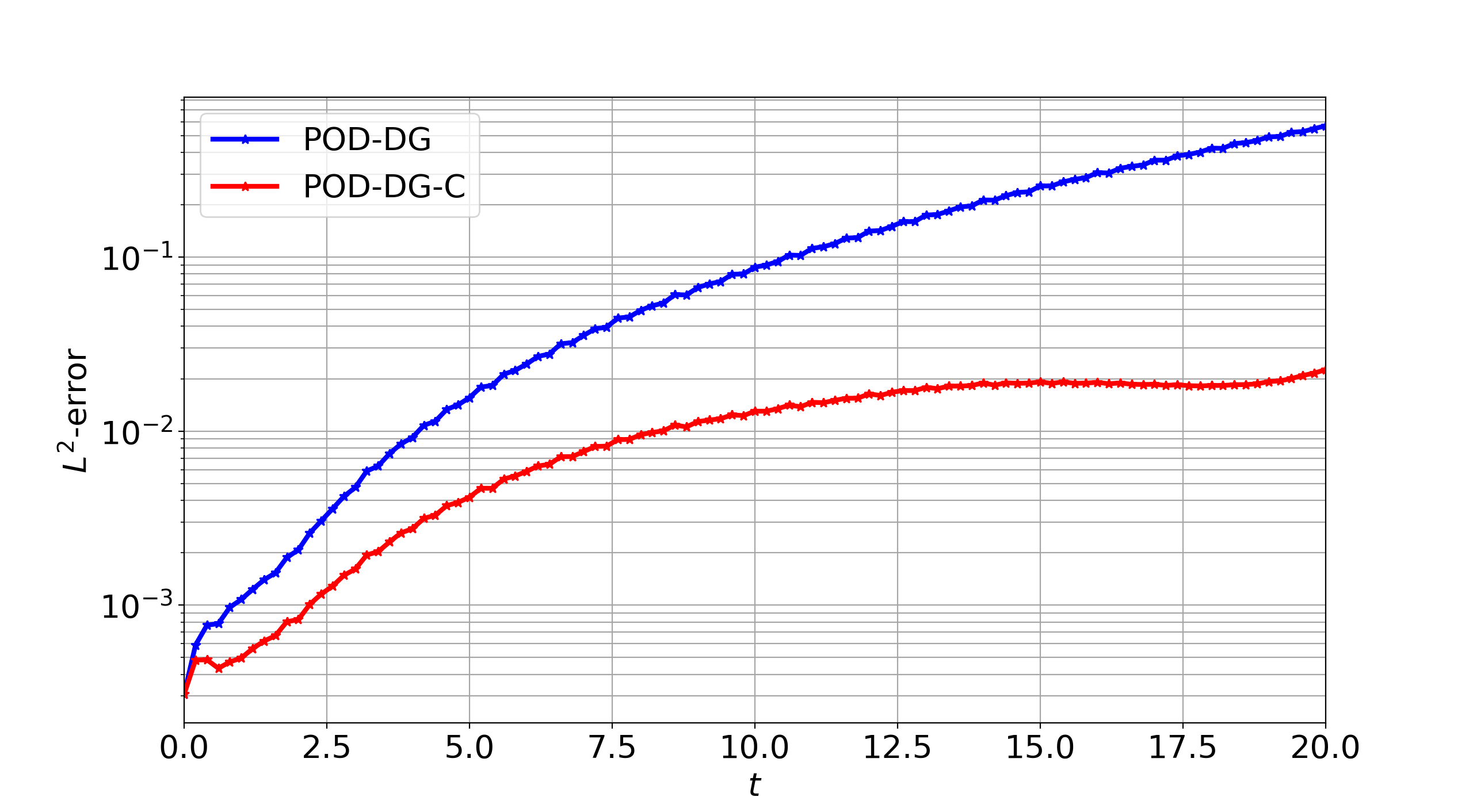}
  \caption{Example 4: time evolution of 
  the $L^2$-velocity error 
$\|\bld
u_h^{fom}-\bld u_h^{rom}\|$.} 
 \label{fig:ns-err2}
\end{figure}

To build the POD bases, we 
collect $501$ snapshots in the time interval
$[0,2]$ taken at equidistant time instance. 
To construct the POD model, we use 10 POD bases which capture about 
99.90\% of the total energy and run the simulation up to time 
$T=20$. We consider the plain POD-DG model and the POD-DG-C model with
$c_1=12$, which is tuned to yield relatively 
smallest $L^2$-error between FOM and ROM solutions at final time. 
Again, we find that adding extra diffusive stabilization in \eqref{closure} does not
improve the results. Hence, results for the POD-DG-CD model will not be
shown.
The time evolution of the $L^2$ velocity error 
$\|\bld
u_h^{fom}-\bld u_h^{rom}\|$ 
is plotted in Figure~\ref{fig:ns-err2}.
We observe that the error for the POD-DG-C model is again an order of
magnitude smaller than that for the plain POD-DG model at time $t=20$.

\begin{figure}[ht!]
 \includegraphics[width=.75\textwidth]{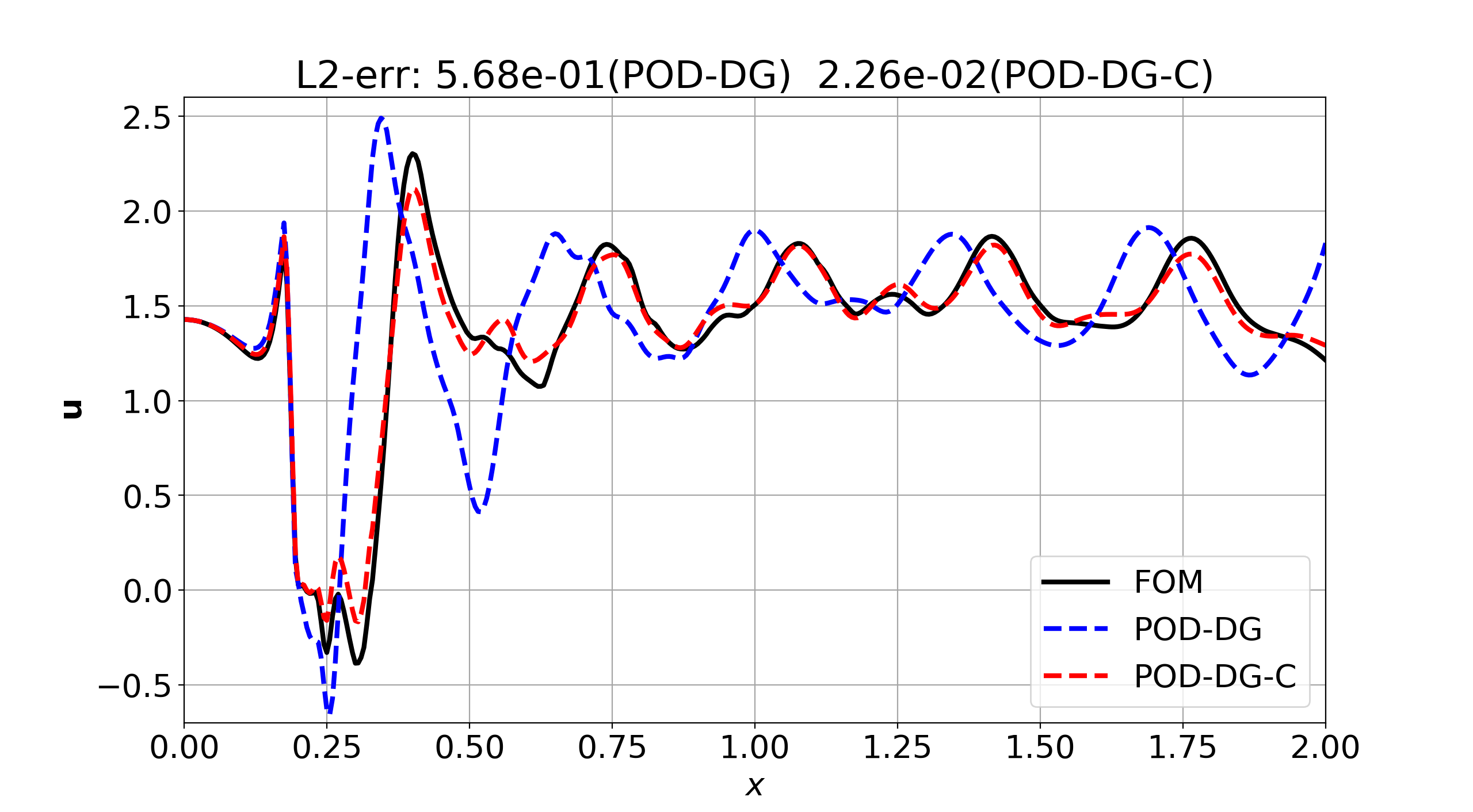}
  \caption{Example 4: x-component of velocity field along cut line 
  $y=0.25$ at time $t=20$.} 
 \label{fig:ns-u21}
\end{figure}

\begin{figure}[ht!]
 \includegraphics[width=.45\textwidth]{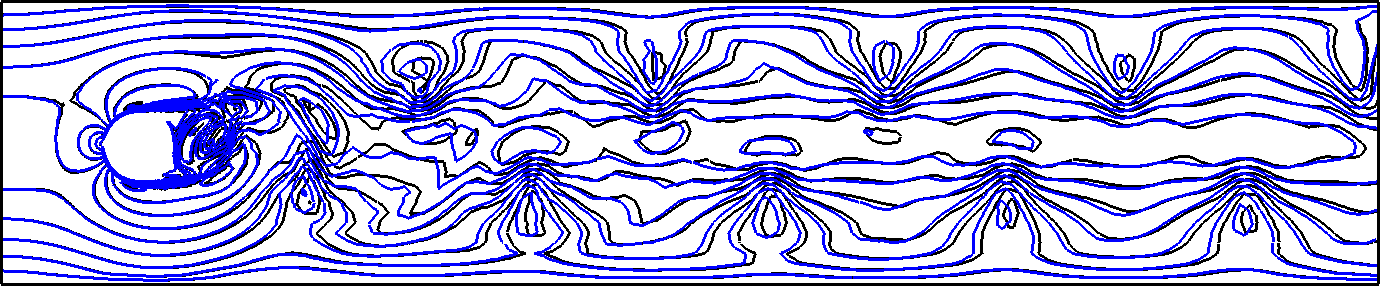}
 \includegraphics[width=.45\textwidth]{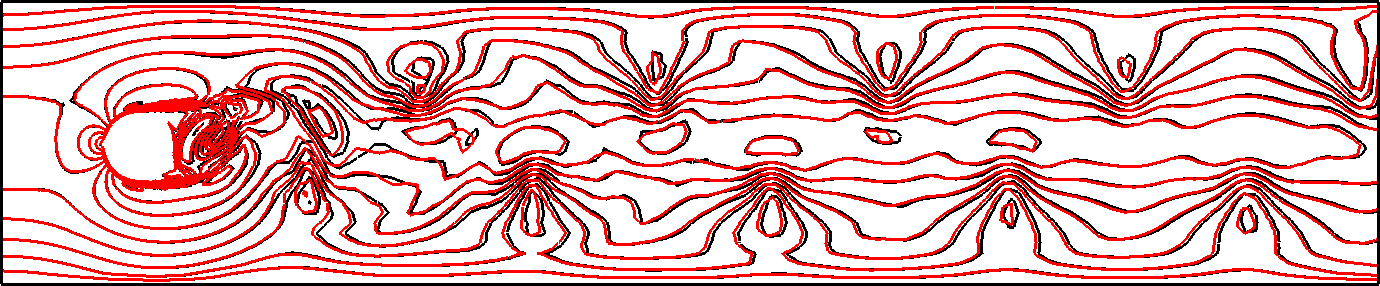}\\[.5ex]
 \includegraphics[width=.45\textwidth]{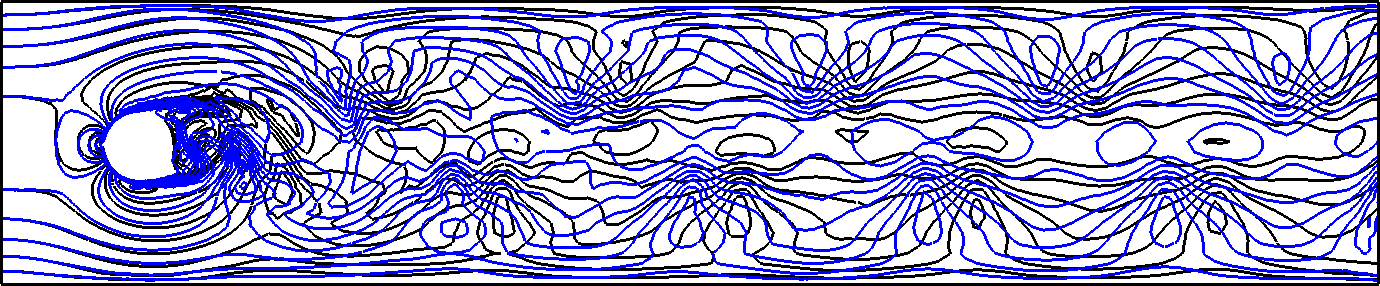}
 \includegraphics[width=.45\textwidth]{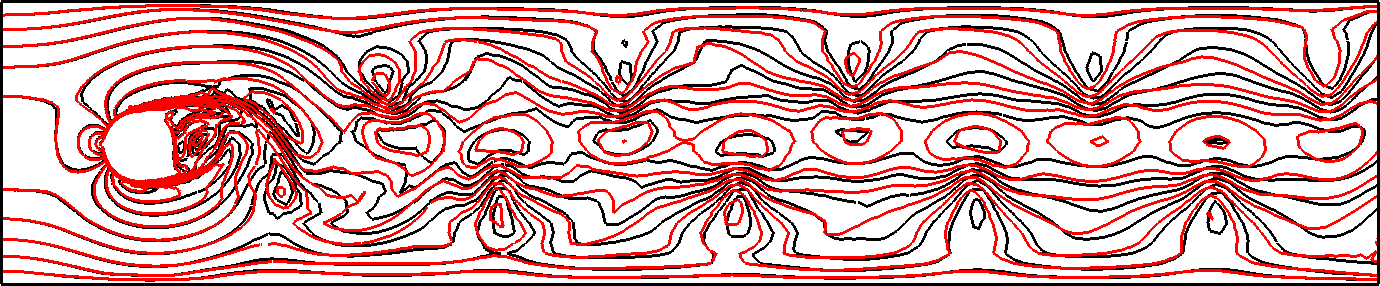}
 \caption{Example 4 (ROMs are 10-dimensional): velocity contour at time $t=5$ (top) and $t=20$
   (bottom). 10 equispaced contour lines from $0$ to $2.4$.
   Black: FOM. Bule: POD-DG. Red: POD-DG-C with $c_1=12$. 
} 
 \label{fig:ns-u2}
\end{figure}

We plot 
the x-component of the velocity field along the 
cut line $y=0.25$ at time $t=20$ in Figure~\ref{fig:ns-u21}, 
and velocity contour lines at time $t=5$ and $t=20$ in
Figure~\ref{fig:ns-u2}. Similar results as those in Example 3 is observed.
In particular, while both models produces similar results at time $t=5$. 
Significant improvement from POD-DG-C model over the plain POD-DG model is
observed for the velocity magnitude contour lines at time $t=20$. This indicates our
POD-DG-C model is more accurate than POD-DG model 
for long time simulations.

\subsection{Example 5. Incompressible Euler: 
double shear layer problem}
In our last example, we consider 
the classical double shear layer problem \cite{BellColella89}.  
We solve the Euler equation \eqref{eq:ns} 
with $\nu=0$ on a periodic  domain $[0,2\pi]\times [0,2\pi]$ 
with an initial condition:
\begin{align}
 u_1(x,y,0) = &\;\left\{
 \begin{tabular}{ll}
$\mathrm{tanh}((y-\pi/2)/\rho)$  & $y\le \pi$\\[1ex]
$\mathrm{tanh}((3\pi/2-y)/\rho)$  & $y> \pi$\\
  \end{tabular}
 \right.,\\
u_2(x,y,0) = &\; \delta \sin(x), 
\end{align}
with $\rho = \pi/15$ and $\delta = 0.05$.

\begin{figure}[ht!]
 \includegraphics[width=.75\textwidth]{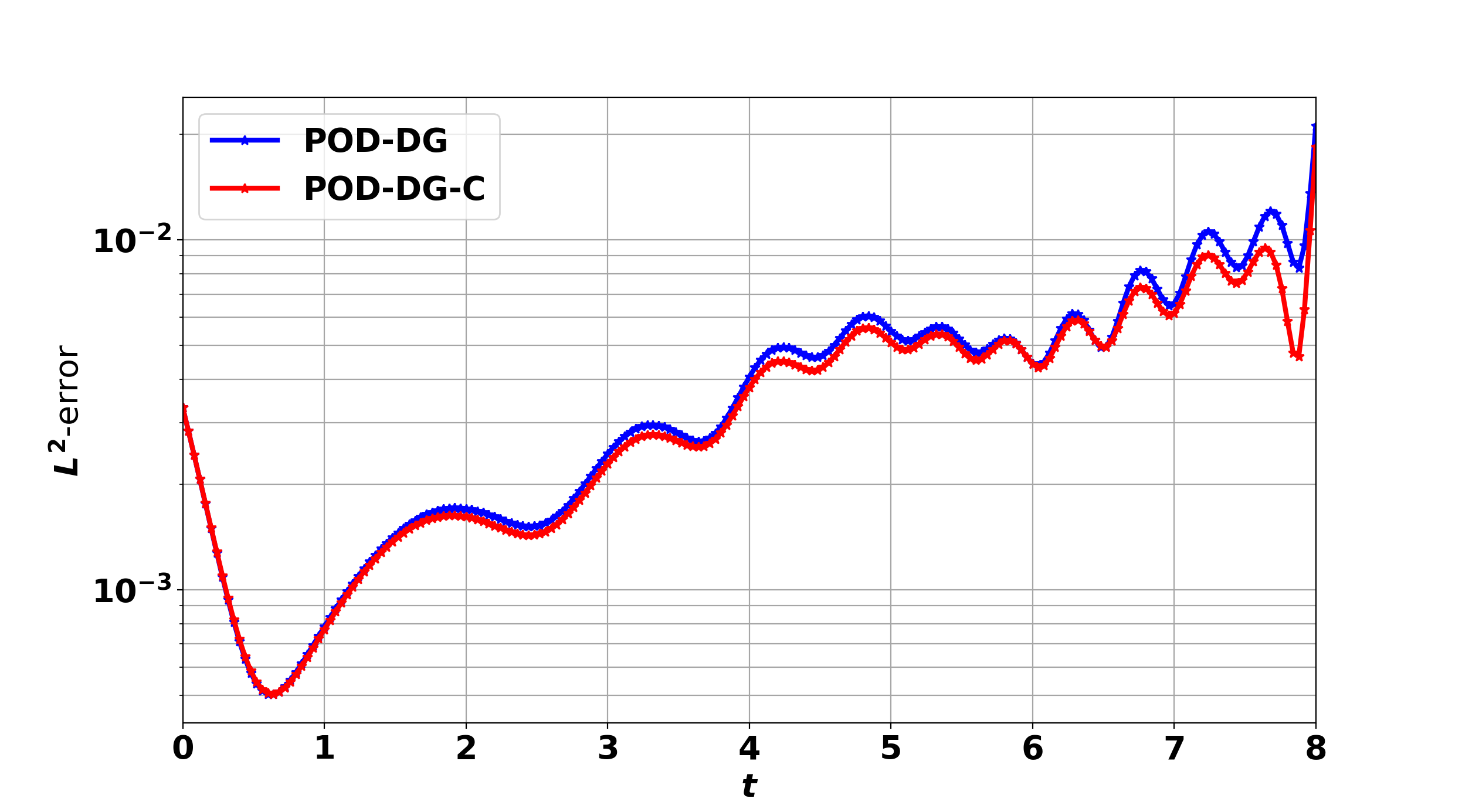}
  \caption{Example 5: time evolution of 
  the $L^2$-velocity error 
$\|\bld
u_h^{fom}-\bld u_h^{rom}\|$.} 
 \label{fig:euler-err}
\end{figure}

For the FOM \eqref{full-n}, we use $P^3$ approximation on fixed uniform
structured triangular meshes with mesh size $2\pi/64$ and run the simulation  up to time $t=8$
with time step size $\Delta t = 0.001$.
To build the POD bases, we collect 401 snapshots in the time interval [0, 8] taken at equidistant
time instance. To build the POD model, we use 10 POD bases which captures
about 99.95\% of the
total energy and run the simulation up to time $T = 8$. We consider the plain POD-DG model and
the POD-DG-C model. The parameter $c_1$ in the POD-DG-C model is tuned to
be $c_1=40$. 
The time evolution of the $L^2$ velocity error 
$\|\bld
u_h^{fom}-\bld u_h^{rom}\|$ 
is plotted in Figure~\ref{fig:euler-err}.
In contract to Examples 3-4, we observe that  the error for both models are
very similar,
which indicates that our current convective stabilization
approach is not effective for the current problem.
We further remark that we also observe no accuracy improvement by
considering the POD-DG-CD model.

\begin{figure}[ht!]
 \includegraphics[width=.45\textwidth]{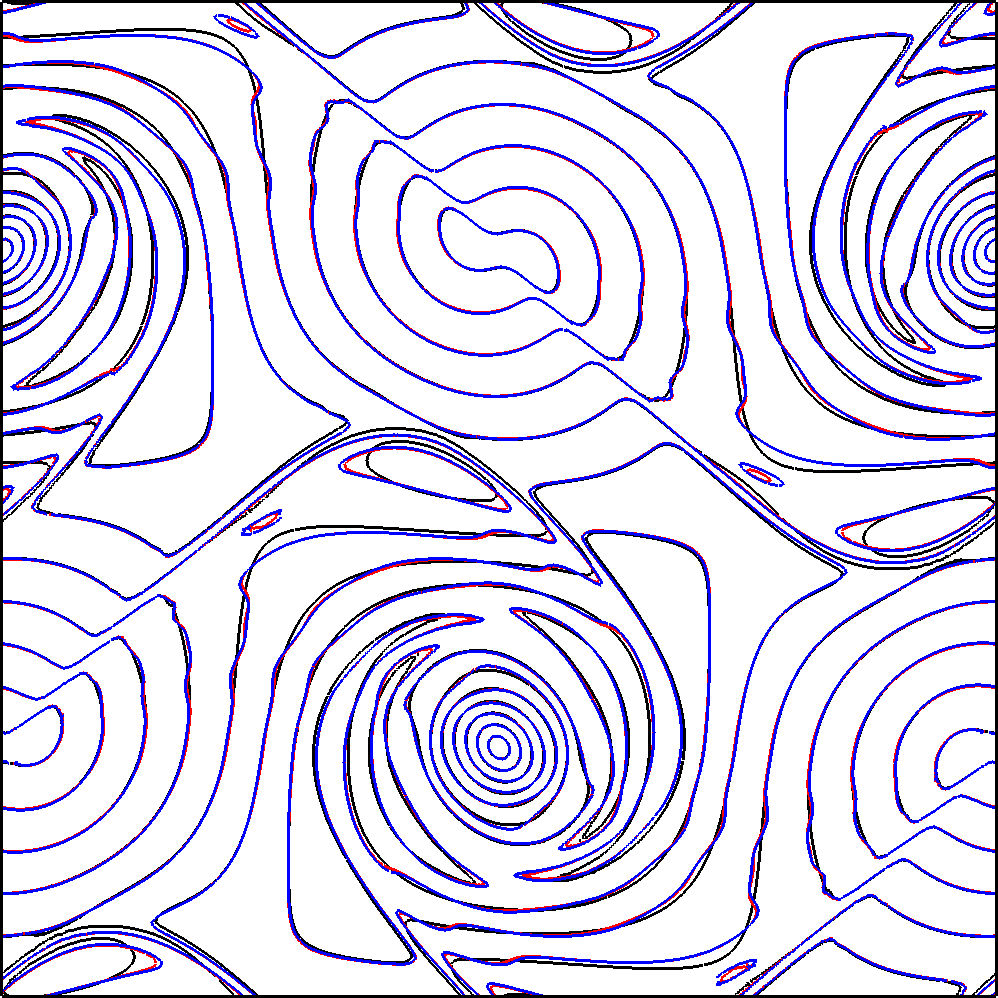}
 \includegraphics[width=.45\textwidth]{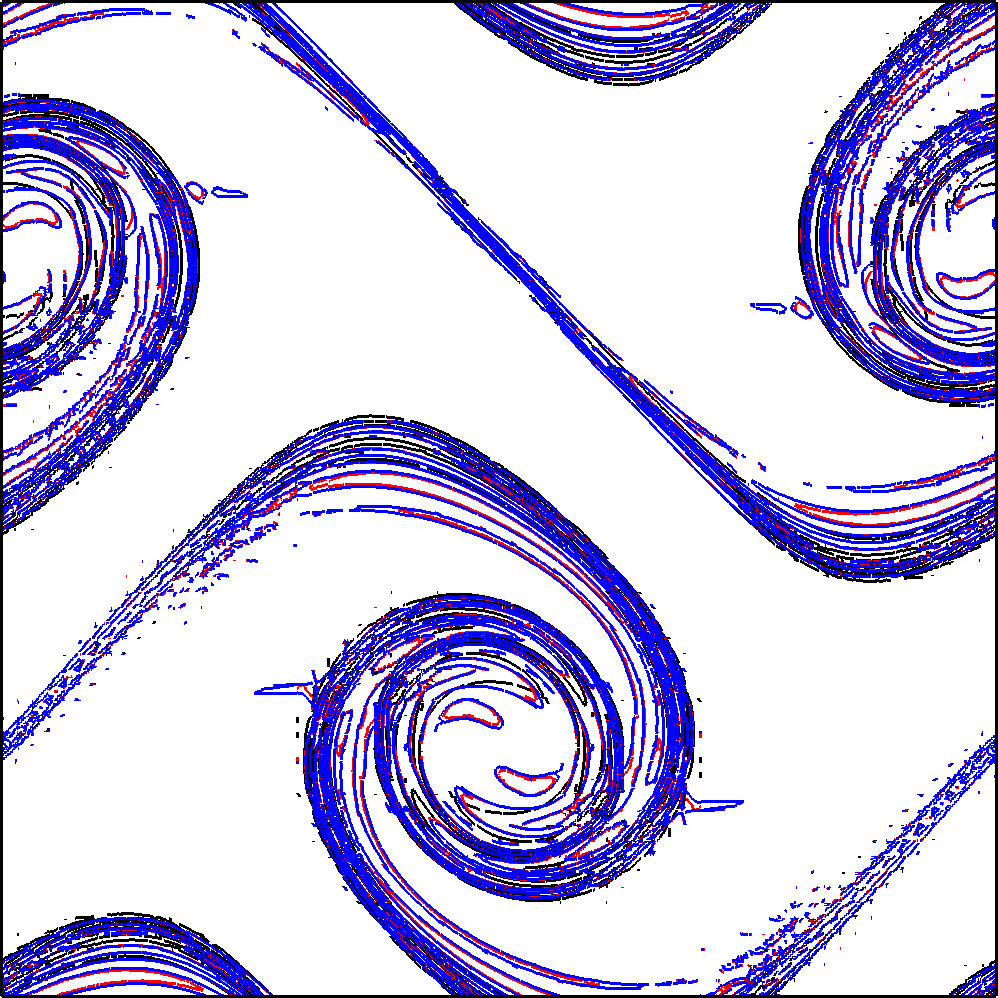}
 \caption{Example 5 (ROMs are 10-dimensional):
   Left:
10 equispaced velocity magnitude contour lines from $0$ to $1.5$. Right: 
10 equispaced vorticity contour lines from $-4.9$ to $4.9$. 
   Black: FOM. Bule: POD-DG. Red: POD-DG-C with $c_1=40$. 
} 
 \label{fig:euler}
\end{figure}

\begin{figure}[ht!]
 \includegraphics[width=.45\textwidth]{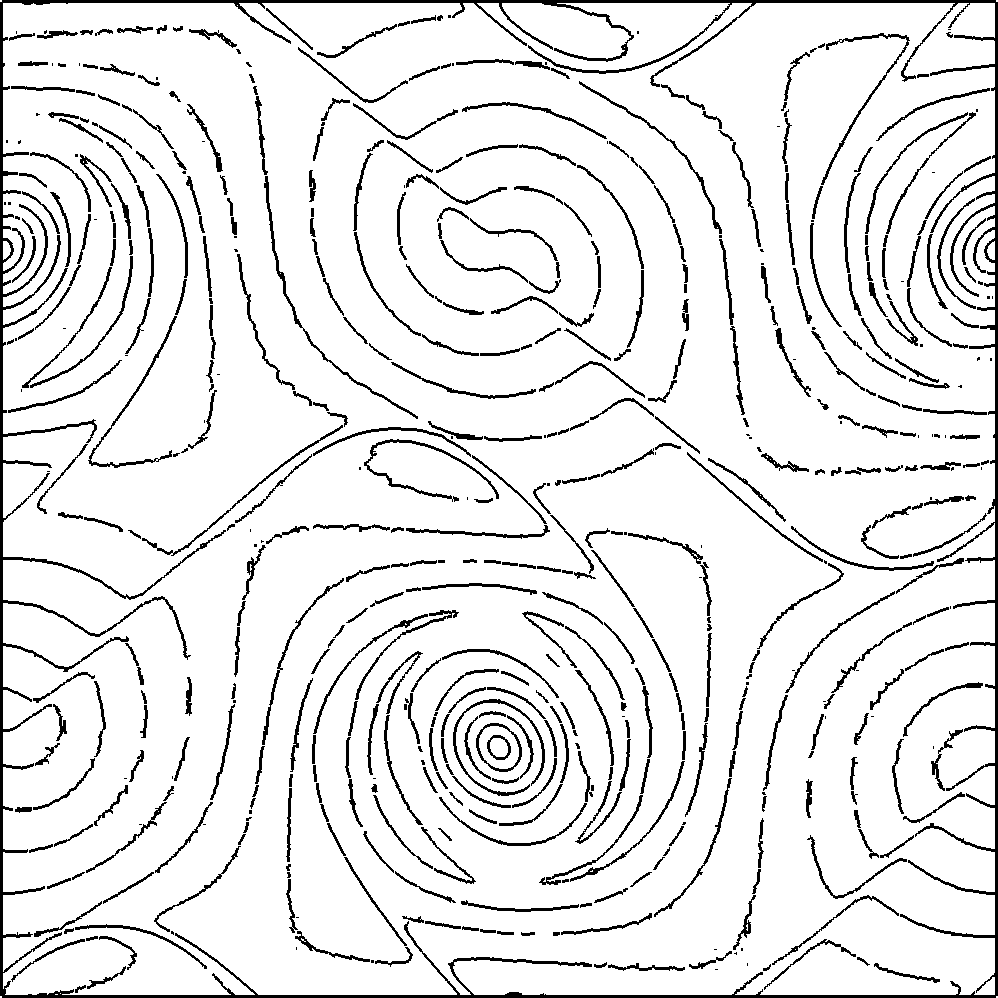}
 \includegraphics[width=.45\textwidth]{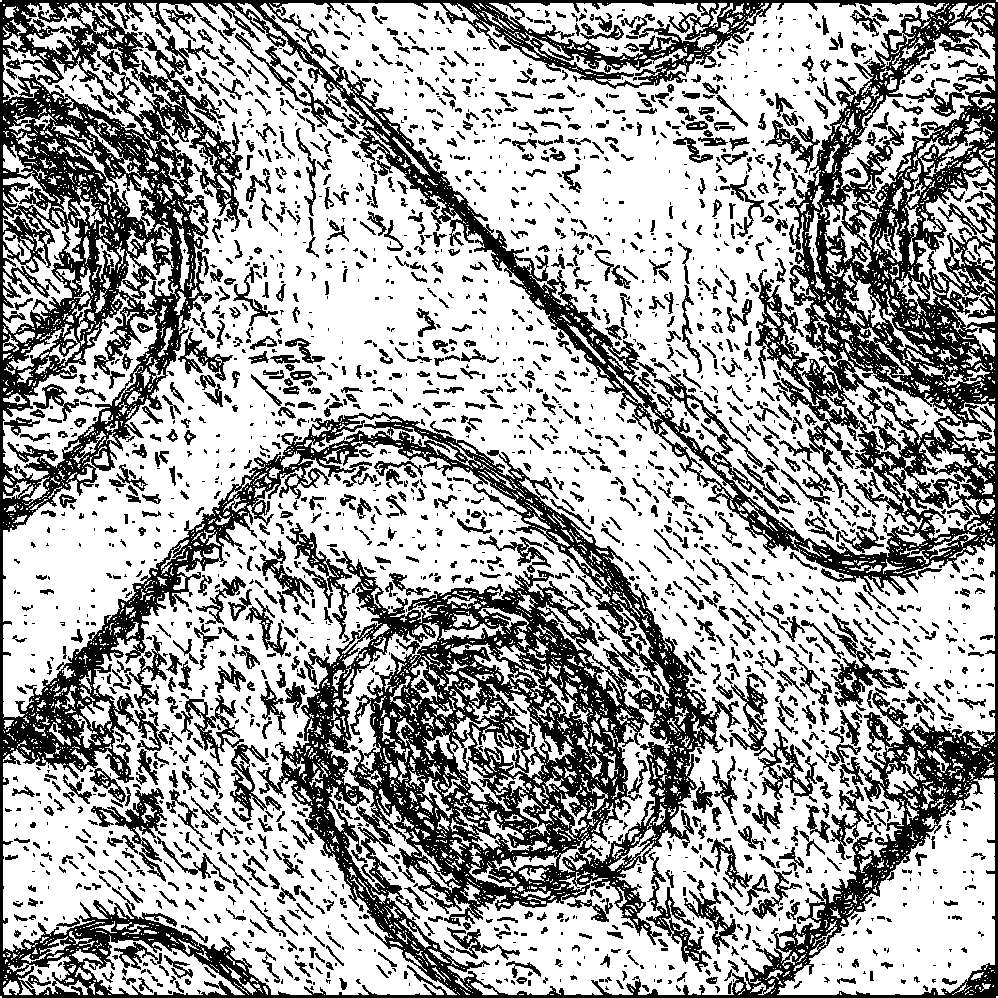}
 \caption{Example 5:
   Left:
10 equispaced velocity magnitude contour lines from $0$ to $1.5$. Right: 
10 equispaced vorticity contour lines from $-4.9$ to $4.9$. 
Results for (conservative) C-FOM.
} 
 \label{fig:ns-X}
\end{figure}

Finally, we plot 
velocity magnitude ($\|\bld u_h\|$) and 
vorticity ($\nabla\times \bld u_h$)
contour lines for 
the two models along with the results for the FOM at final time $t=8$
in
Figure~\ref{fig:euler}.  
It can be observed that the results for both models are very similar and
are close to the FOM results. 
This is a rather surprising result as the POD-DG model does not introduce
any spatial numerical dissipation, yet its vorticity approximation is still
free from large oscillations. For comparison, we also plot in
Figure~\ref{fig:ns-X} the vorticity
approximations for the FOM \eqref{full-n} with 
the upwinding convection operator $\bld {\mathcal{C}}_h^{dg}$ replaced by
the conservative version 
$\bld {\widetilde{\mathcal{C}}}_h^{dg}$ in \eqref{conv-c}, which we denote
as
C-FOM.
It is clear that the (conservative) POD-DG model has better stability
property than C-FOM. We conjecture the reason for the superior performance
of POD-DG model over C-FOM is that the global POD bases obtained from FOM \eqref{full-n}
might have some extra built-in stabilization properties.

\section{Conclusion}
\label{sec:conclude}
We have presented a POD-DG reduced order model 
for the viscous Burgers' equation and
the incompressible Navier-Stokes equations discretized using an IMEX HDG/DG
scheme. A novel offline-computable closure model was introduced for the POD-DG ROM
which further improves its stability and accuracy. 
Numerical results show the superior performance of the proposed closure
model comparing with a plain POD-DG scheme without the closure model.
In future work, we will pursue in the same direction and investigate the proposed model in the parametrized flow problems 
with applications in flow control and optimization. 

\bibliographystyle{siam}

\begin{thebibliography}{10}

\bibitem{antonietti2016discontinuous}
{\sc P.~F. Antonietti, P.~Pacciarini, and A.~Quarteroni}, {\em A discontinuous
  {G}alerkin reduced basis element method for elliptic problems}, ESAIM:
  Mathematical Modelling and Numerical Analysis, 50 (2016), pp.~337--360.

\bibitem{antoulas2020interpolatory}
{\sc A.~Antoulas, C.~Beattie, and S.~G{\"u}{\u{g}}ercin}, {\em Interpolatory
  methods for model reduction}, 2020.

\bibitem{antoulas2005approximation}
{\sc A.~C. Antoulas}, {\em Approximation of large-scale dynamical systems},
  vol.~6, Siam, 2005.

\bibitem{AscherRuuthWetton93}
{\sc U.~M. Ascher, S.~J. Ruuth, and B.~T.~R. Wetton}, {\em Implicit-explicit
  methods for time-dependent partial differential equations}, SIAM J. Numer.
  Anal., 32 (1995), pp.~797--823.

\bibitem{ballarin2015supremizer}
{\sc F.~Ballarin, A.~Manzoni, A.~Quarteroni, and G.~Rozza}, {\em Supremizer
  stabilization of {POD}--{G}alerkin approximation of parametrized steady
  incompressible {N}avier--{S}tokes equations}, International Journal for
  Numerical Methods in Engineering, 102 (2015), pp.~1136--1161.

\bibitem{BellColella89}
{\sc J.~B. Bell, P.~Colella, and H.~M. Glaz}, {\em A second-order projection
  method for the incompressible {N}avier-{S}tokes equations}, J. Comput. Phys.,
  85 (1989), pp.~257--283.

\bibitem{benner2015survey}
{\sc P.~Benner, S.~Gugercin, and K.~Willcox}, {\em A survey of projection-based
  model reduction methods for parametric dynamical systems}, SIAM review, 57
  (2015), pp.~483--531.

\bibitem{caiazzo2014numerical}
{\sc A.~Caiazzo, T.~Iliescu, V.~John, and S.~Schyschlowa}, {\em A numerical
  investigation of velocity--pressure reduced order models for incompressible
  flows}, Journal of Computational Physics, 259 (2014), pp.~598--616.

\bibitem{carlberg2013gnat}
{\sc K.~Carlberg, C.~Farhat, J.~Cortial, and D.~Amsallem}, {\em The {GNAT}
  method for nonlinear model reduction: effective implementation and
  application to computational fluid dynamics and turbulent flows}, Journal of
  Computational Physics, 242 (2013), pp.~623--647.

\bibitem{chaturantabut2010nonlinear}
{\sc S.~Chaturantabut and D.~C. Sorensen}, {\em Nonlinear model reduction via
  discrete empirical interpolation}, SIAM Journal on Scientific Computing, 32
  (2010), pp.~2737--2764.

\bibitem{Cockburn16}
{\sc B.~Cockburn}, {\em Static condensation, hybridization, and the devising of
  the {HDG} methods}, in Building bridges: connections and challenges in modern
  approaches to numerical partial differential equations, vol.~114 of Lect.
  Notes Comput. Sci. Eng., Springer, [Cham], 2016, pp.~129--177.

\bibitem{gunzburger2003perspectives}
{\sc M.~D. Gunzburger}, {\em Perspectives in flow control and optimization},
  vol.~5, Siam, 2003.

\bibitem{hesthaven2016certified}
{\sc J.~S. Hesthaven, G.~Rozza, B.~Stamm, et~al.}, {\em Certified reduced basis
  methods for parametrized partial differential equations}, vol.~590, Springer,
  2016.

\bibitem{Lehrenfeld10}
{\sc C.~Lehrenfeld}, {\em Hybrid {Discontinuous} {Galerkin} methods for solving
  incompressible flow problems}.
\newblock Diploma Thesis, MathCCES/IGPM, RWTH Aachen, 2010.

\bibitem{LehrenfeldSchoberl16}
{\sc C.~Lehrenfeld and J.~Sch\"{o}berl}, {\em High order exactly
  divergence-free hybrid discontinuous {G}alerkin methods for unsteady
  incompressible flows}, Comput. Methods Appl. Mech. Engrg., 307 (2016),
  pp.~339--361.

\bibitem{peherstorfer2016data}
{\sc B.~Peherstorfer and K.~Willcox}, {\em Data-driven operator inference for
  nonintrusive projection-based model reduction}, Computer Methods in Applied
  Mechanics and Engineering, 306 (2016), pp.~196--215.

\bibitem{quarteroni2015reduced}
{\sc A.~Quarteroni, A.~Manzoni, and F.~Negri}, {\em Reduced basis methods for
  partial differential equations: an introduction}, vol.~92, Springer, 2015.

\bibitem{SanIliescu13}
{\sc O.~San and T.~Iliescu}, {\em Proper orthogonal decomposition closure
  models for fluid flows: {B}urgers equation}, Int. J. Numer. Anal. Model. Ser.
  B, 5 (2014), pp.~217--237.

\bibitem{LD96}
{\sc M.~Sch\"afer, S.~Turek, F.~Durst, K.~E., and R.~R.}, {\em Benchmark
  computations of laminar flow around a cylinder}, Flow simulation with
  high-performance computers II 1996; :547–566.

\bibitem{Schoberl16}
{\sc J.~Sch{\"o}berl}, {\em {C}++11 {I}mplementation of {F}inite {E}lements in
  {NGS}olve}, 2014.
\newblock {ASC Report 30/2014, Institute for Analysis and Scientific Computing,
  Vienna University of Technology}.

\bibitem{shen2019hdg}
{\sc J.~Shen, J.~R. Singler, and Y.~Zhang}, {\em {HDG}--{POD} reduced order
  model of the heat equation}, Journal of Computational and Applied
  Mathematics, 362 (2019), pp.~663--679.

\bibitem{Sirovich87}
{\sc L.~Sirovich}, {\em Turbulence and the dynamics of coherent structures.
  {I}. {C}oherent structures}, Quart. Appl. Math., 45 (1987), pp.~561--571.

\bibitem{uzunca2017energy}
{\sc M.~Uzunca and B.~Karas{\"o}zen}, {\em Energy stable model order reduction
  for the {A}llen-{C}ahn equation}, in Model Reduction of Parametrized Systems,
  Springer, 2017, pp.~403--419.

\bibitem{wang2015nonlinear}
{\sc Z.~Wang}, {\em Nonlinear model reduction based on the finite element
  method with interpolated coefficients: semilinear parabolic equations},
  Numerical Methods for Partial Differential Equations, 31 (2015),
  pp.~1713--1741.

\bibitem{yano2019discontinuous}
{\sc M.~Yano}, {\em Discontinuous {G}alerkin reduced basis empirical quadrature
  procedure for model reduction of parametrized nonlinear conservation laws},
  Advances in Computational Mathematics,  (2019), pp.~1--34.

\end{thebibliography}

\end{document}